\documentclass[12pt]{article}
\usepackage{amsmath,amssymb,color}

\definecolor{mycolorred}{rgb}{1, 0, 0}

\def\P{{\mathbb P}}
\def\R{{\mathbb R}}
\def\E{{\mathbb E}}

\def\N{{\mathbb N}}

\def\<{\langle}
\def\>{\rangle}

\newtheorem{theorem}{Theorem}[section]

\newtheorem{corollary}[theorem]{Corollary}

\newtheorem{lemma}[theorem]{Lemma}

\newtheorem{proposition}[theorem]{Proposition}
\newtheorem{remark}[theorem]{Remark}

\numberwithin{equation}{section}

\begin{document}

\title{Regularity of probability laws by using an interpolation method}
\author{ \textsc{Vlad Bally}\thanks{%
Laboratoire d'Analyse et de Math\'ematiques Appliqu\'ees, UMR 8050,
Universit\'e Paris-Est Marne-la-Vall\'ee, 5 Bld Descartes, Champs-sur-Marne,
77454 Marne-la-Vall\'ee Cedex 2, France. Email: \texttt{bally@univ-mlv.fr} }%
\smallskip \\
\textsc{Lucia Caramellino}\thanks{%
Dipartimento di Matematica, Universit\`a di Roma - Tor Vergata, Via della
Ricerca Scientifica 1, I-00133 Roma, Italy. Email: \texttt{%
caramell@mat.uniroma2.it}}\smallskip\\
}
\maketitle

\parindent 0pt

{\textbf{Abstract.}}
We study the problem of the existence and regularity of a probability density in an abstract framework based on a ``balancing'' with approximating absolutely continuous laws.  Typically, the absolutely continuous property for the approximating laws can be proved by standard techniques from Malliavin calculus whereas for the law of interest no Malliavin integration by parts formulas are available. Our results are strongly based on the use of suitable Hermite polynomial series expansions and can be merged into the theory of interpolation spaces. We then apply the results  to the solution to a stochastic differential equation with a local H\"ormander condition or to the solution to the stochastic heat equation, in both cases under weak conditions on the coefficients relaxing the standard Lipschitz or H\"older continuity requests.

\medskip


\medskip


\tableofcontents

\section{Introduction}

P. Malliavin has built a stochastic differential calculus which allows one
to prove integration by parts formulas of type $\E(\partial _{\alpha }\phi
(F))=\E(\phi (F)H_{\alpha }(F))$ and used them in order to study the
regularity of the law of $F.$ Here $F$ is a functional on the Wiener space.
Roughly speaking the strategy is the following: on takes a sequence of
simple functionals $F_{n}\rightarrow F,$ defines the differential operators $%
DF_{n}$ and then, if $F\in DomD,$ defines $DF=\lim_{n}DF_{n}.$ And this
infinite dimensional differential calculus allows one to prove the
integration by parts formulas. But one may proceed in a different way: using
the finite dimensional calculus associated to simple functionals one proves $%
\E(\partial _{\alpha }\phi (F_{n}))=\E(\phi (F_{n})H_{\alpha }(F_{n})),$ use
it in order to get estimates of the Fourier transform and then pass to the
limit. Of course, if everything works well when passing to the limit, this
is more or less the same. But the interesting point is that one may use this
strategy even if $\E(\phi (F_{n})H_{\alpha }(F_{n}))$ does not converge, so
for $F\notin DomD.$ In fact, consider a random variable $F$ and a sequence
of functionals $F_{n},n\in \N$ such that $\E\left\vert F-F_{n}\right\vert
\rightarrow 0$ and suppose that some integration by parts formulas hold for
each $F_{n}$ (does not matter how one obtains them). Then one can proceed as
follows (for simplicity we consider the one dimensional case). Let $\widehat{%
p}_{F_{n}}(\xi )=\E(e^{i\xi F_{n}})$ be the Fourier transform of $F_{n}.$
Since $\partial _{x}^{k}e^{i\xi F_{n}}=(i\xi )^{k}e^{i\xi F_{n}}$ one may
use $k$ integration by parts and obtains
\begin{equation*}
\widehat{p}_{F_{n}}(\xi )=\frac{1}{(i\xi )^{k}}\E(\partial _{x}^{k}e^{i\xi
F_{n}})=\frac{1}{(i\xi )^{k}}\E(e^{i\xi F_{n}}H_{k}(F_{n})).
\end{equation*}%
Then one writes
\begin{eqnarray}
\left\vert \widehat{p}_{F}(\xi )\right\vert &=&\left\vert \widehat{p}%
_{F}(\xi )-\widehat{p}_{F_{n}}(\xi )\right\vert +\frac{1}{\left\vert \xi
\right\vert ^{k}}\left\vert \E(e^{i\xi F_{n}}H_{k}(F_{n}))\right\vert  \notag
\\
&\leq &\left\vert \xi \right\vert \E\left\vert F-F_{n}\right\vert +\frac{1}{%
\left\vert \xi \right\vert ^{k}}\E\left\vert H_{k}(F_{n})\right\vert .
\label{I1}
\end{eqnarray}%
So if one succeeds to get a good balance between $\E\left\vert
F-F_{n}\right\vert \downarrow 0$ and $\E\left\vert H_{k}(F_{n})\right\vert
\uparrow \infty $ one may obtain good estimates of $\left\vert \widehat{p}%
_{F}(\xi )\right\vert $ and this implies the regularity of the law of $F.$

This argument originates from \cite{bib:[F]} and has been used in several
recent papers: see \cite{bib:[BCl1]}, \cite{bib:[BCl2]}, \cite{bib:[BF]},
\cite{bib:[De]} and \cite{bib:[FP]}. Notice however that this method depends
on the dimension: the weaker condition which gives a density of the law of $%
F $ is $\int_{\R^{d}}\left\vert \widehat{p}_{X_{T}}(\xi )\right\vert
^{2}dx<\infty $ and of course, this depends on $d.$ Let us give a simple but
significant example (see \cite{bib:[FP]}). Consider a $d$ dimensional
diffusion process $dX_{t}=\sum_{j=1}^{N}\sigma
_{j}(X_{t})dW_{t}^{j}+b(X_{t})dt.$ In \cite{bib:[FP]} one assumes that the
coefficients $\sigma _{j}$ and $b$ are H\"{o}lder continuous of order $h>0$
and tries to prove that the law of $X_{T}$ is absolutely continuous. One
takes $\delta >0$ and defines $X_{t}^{\delta }$ to be equal to $X_{t}$ for $%
t\leq T-\delta $ and $X_{t}^{\delta }=X_{T-\delta }+\sum_{j=1}^{N}\sigma
_{j}(X_{T-\delta })(W_{t}^{j}-W_{T-\delta }^{j})$ for $t\geq T-\delta .$ It
is easy to see that $\E\left\vert X_{T}-X_{T}^{\delta }\right\vert \leq
C\delta ^{\frac{1}{2}(1+h)}.$ On the other hand if $\sigma \sigma ^{\ast
}\geq c$ then conditionally to $X_{T-\delta },$ $X_{T}^{\delta }$ is a non
degenerated Gaussian random variable. Using elementary integration by parts
one may prove that for every multi index $\alpha $ of length $k,$ $%
\E(\partial _{\alpha }\phi (X_{T}^{\delta }))=\E(\phi (X_{T}^{\delta
})H_{\alpha ,\delta })$ and $\E\left\vert H_{\alpha ,\delta }\right\vert \leq
C\delta ^{-k/2}.$ Using (\ref{I1})
\begin{equation*}
\left\vert \widehat{p}_{X_{T}}(\xi )\right\vert \leq \left\vert \xi
\right\vert \E\left\vert X_{T}-X_{T}^{\delta }\right\vert +\frac{1}{%
\left\vert \xi \right\vert ^{k}}\E\left\vert H_{\alpha ,\delta }\right\vert
\leq C\left\vert \xi \right\vert \delta ^{\frac{1}{2}(1+h)}+\frac{C}{%
\left\vert \xi \right\vert ^{k}\delta ^{k/2}}.
\end{equation*}%
We fix $\xi $ and we choose $\delta =\left\vert \xi \right\vert
^{-2(k+1)/(k+1+h)}$ in order to optimize the above relation. Then
\begin{equation*}
\left\vert \widehat{p}_{X_{T}}(\xi )\right\vert \leq C\left\vert \xi
\right\vert ^{-\frac{hk}{1+h+k}}.
\end{equation*}%
If $d=1$ and $h>\frac{1}{2}$ one may choose $k$ sufficiently large in order
to get $\frac{hk}{1+h+k}>\frac{1}{2}$ and so $\int_{\R^{d}}\left\vert
\widehat{p}_{X_{T}}(\xi )\right\vert ^{2}dx<\infty .$ But for $d\geq 2,$
even if $h=1, $ one fails to prove that the above integral is finite. So
this approach is successful just for $d=1$ and $h>\frac{1}{2}.$

\smallskip

The aim of this paper is to obtain a more performing balance which
essentially does not depend on the dimension. The main results are Theorem %
\ref{2E} in Section \ref{sect-results} and Theorem \ref{2F} in Section \ref%
{sect-IBP}. Our estimates are based on a development in Hermite series
(instead of the Fourier transform) presented in Section \ref{sect-Hermite}.
In this framework we use a powerful result concerning the regularity of a
mixture of Hermite functions - see Theorem \ref{5} in Section \ref%
{sect-Hermite}. This result has been proved in \cite{bib:[E]} in the one
dimensional case then in \cite{bib:[Dz]} for the multi dimensional case. We
use a variant given in \cite{bib:[PY]} Corollary 2.3.

\smallskip

As an application of Theorem \ref{2E} \ we are able to improve the above
mentioned result of \cite{bib:[FP]} in the following way. For an open domain
$D\subset \R^{d}$ we denote by $C_{\log }(D)$ the class of functions $%
f:D\rightarrow \R$ for which there exists $C,h>0$ such that $\left\vert
f(x)-f(y)\right\vert \leq C\left\vert \ln \left\vert x-y\right\vert
\right\vert ^{-h}$ for every $x,y\in D.$ We fix $y_{0}\in \R^{d}$ and $r>0$ \
and we suppose that the coefficients of the SDE presented above verify $%
\sigma _{j}\in C_{\log }(B_{r}(y_{0})),j=1,...,N$ (the drift coefficient $b$
is just measurable with linear growth). We also assume that $\sigma \sigma
^{\ast }(y_{0})>0.$ Finally we consider an open domain $\Gamma \subset \R^{d}$
and we denote by $\tau $ the exit time from $\Gamma .$ Our result (see
Theorem \ref{9}) says that for $y_{0}\in \Gamma $ which verifies the above
hypothesis the law of $X_{T\wedge \tau }$ is absolutely continuous with
respect to the Lebesgue measure in a neighborhood of $y_{0}.$ And this is
true for any dimension $d$ and every $h>0.$ Notice that even if the
coefficients are smooth one cannot use directly the Malliavin calculus
because, due to the stopping time $\tau ,$ $X_{T\wedge \tau }$ is not
differentiable in Malliavin sense. In \cite{bib:[C]} it is given a variant
of Malliavin calculus which permits to handle SDE's with boundary
conditions - but there the coefficients are smooth (while here they are just
in $C_{\log }(B_{r}(y_{0})))$. Finally in Section 5 we prove a similar absolute continuity result for
solutions of the stochastic heat equation.

\smallskip

We mention also that recently Debussche and Romito \cite{DR} introduced an
alternative approach, based on Besov space techniques, which enables one to
prove the above result for H\"{o}lder continuous coefficients. This
technique has already been used by Fournier \cite{F1} in order to study the
regularity of the $3$-dimensional Boltzman equation.

\smallskip

After having done the work concerning probability measures we realized that
this fits in the more general theory of interpolation spaces and extends to
distributions (we thank to D. Elworthy who remarked this). In this framework
our criterion ensures that a given distribution (in particular a probability
measure) belongs to a certain interpolation space between a distribution
space and a weighted Sobolev space. And the work done in the paper consists
in proving that such an interpolation space is in between two Sobolev
spaces. This gives the regularity result. There already exist a certain
number of results concerning interpolation between negative Sobolev spaces
and Sobolev spaces. But our result does not seem to fit in this framework.
The reason is that when working with Sobolev spaces one employs $L^{p}$
norms while the distribution space that we consider is defined in terms of $%
L^{\infty }$ norms. And the case of $L^{\infty }$ norms appears as a limit
case in interpolation theory and is more delicate to treat (see Triebel \cite%
{T} or Bennet and Sharepley \cite{bib:[BS]}).

\smallskip

The paper is organized as follows. In Section \ref{sect-results} we state
the main result which is a criterion of regularity for general finite
measures. We prove it in Section \ref{sect-Hermite}, and in Section \ref%
{sect-IBP} we give an alternative regularity criterion using integration by
parts formulas. In Section \ref{sect-interp} we discuss the link with
interpolation spaces. In Section \ref{sect-diffproc} we give two examples.
The first one concerns diffusion processes with coefficients in $C_{\log }(D)$ under
an ellipticity assumption (the example presented above). It turns out that
in this case one does not need to use integration by parts: the analysis
relies on the explicit Gaussian density. The second example concerns
diffusions with more regular coefficients which verify a local H\"{o}rmander
condition. In this case the integration by parts formula from Malliavin
calculus is used. Finally, in Section \ref{sect-heat} we prove a regularity
result for the stochastic heat equation introduced by Walsh \cite{bib:[W]}.
This considerably improves a previous result of Bally and Pardoux \cite%
{bib:[BP]}. In the Appendix we discuss and prove some properties related to
interpolation spaces.

\section{Criterion for the regularity of a probability law}

\subsection{Notation and main results}

\label{sect-results}

We work on $\R^{d}$ and we denote by $\mathcal{M}$ the set of the finite
signed measures on $\R^{d}$ with the Borel $\sigma $ algebra. Moreover $%
\mathcal{M}_{a}\mathcal{\subset M}$ is the set of the measures which are
absolutely continuous with respect to the Lebesgue measure. For $\mu \in
\mathcal{M}_{a}$ we denote by $p_{\mu }$ the density of $\mu $ with respect
to the Lebesgue measure. And for a measure $\mu \in \mathcal{M}$ we denote
by $L_{\mu }^{p}$\ the space of the measurable functions $f:\R^{d}\rightarrow \R$ such that $\int \left\vert f\right\vert ^{p}d\mu <\infty .$ For $f\in
L_{\mu }^{1}$ we denote $f\mu $ the measure $(f\mu )(A)=\int_{A}fd\mu .$ For
a bounded function $\phi :\R^{d}\rightarrow \R$ we denote $\mu \ast \phi $ the
measure defined by $\int fd\mu \ast \phi =\int f\ast \phi d\mu =\int \int
\phi (x-y)f(y)dyd\mu (x).$ Then $\mu \ast \phi \in \mathcal{M}_{a}$ and $%
p_{\mu \ast \phi }(x)=\int \phi (x-y)d\mu (y).$

We denote by $\alpha =(\alpha _{1},...,\alpha _{d})\in \N^{d}$ a multi index
and we put $\left\vert \alpha \right\vert =\sum_{i=1}^{d}\alpha _{i}.$ Here $%
\N=\{0,1,2,...\}$ are the non negative integers and we put $\N_{\ast
}=\N\setminus \{0\}.$ For a multi index $\alpha $ with $\left\vert \alpha
\right\vert =k$ we denote $\partial _{\alpha }$ the corresponding derivative
that is $\partial _{x_{1}}^{\alpha _{1}}...\partial _{x_{d}}^{\alpha _{d}}$
with the convention that $\partial _{x_{i}}^{\alpha _{i}}f=f$ if $\alpha
_{i}=0.$ In particular if $\alpha $ is the null multi index then $\partial
_{\alpha }f=f.$

We denote by $\left\Vert f\right\Vert _{p}=(\int \left\vert f(x)\right\vert
^{p}dx)^{1/p},p\geq 1$ and $\left\Vert f\right\Vert _{\infty }=\sup_{x\in
\R^{d}}\left\vert f(x)\right\vert .$ Then $L^{p}=\{f:\left\Vert f\right\Vert
_{p}<\infty \}$ are the standard $L^{p}$ spaces with respect to the Lebesgue
measure.

In the following we will work in Orlicz spaces so we briefly recall the
notation and the results we will use (we refer to \cite{bib:[K.M-W]}). A
function $e:\R\rightarrow \R_{+}$ is a Young function if it is symmetric,
strictly convex, non negative and $e(0)=0.$ In the following we will
consider a Young function which has the two supplementary properties:%
\begin{align}
i)\quad & \text{there exists }\lambda >0\text{ such that }e(2s)\leq \lambda
e(s),  \label{Y} \\
ii)\quad & s \mapsto \frac{e(s)}{s}\text{ is non decreasing.}  \notag
\end{align}


The property $i)$ is known as the $\Delta _{2}$ condition or doubling
condition (see \cite{bib:[K.M-W]}). Through the whole paper we work with
Young functions which satisfy (\ref{Y}). We denote by $\mathcal{E}$ the
space of these functions.

For $e\in\mathcal{E}$ and $f:\R^{d}\rightarrow \R$, we define the norm
\begin{equation}
\left\Vert f\right\Vert _{(e)}=\inf \{c>0:\int e(\frac{1}{c}f(x))dx\leq 1\}.
\label{O1}
\end{equation}%
This is the so called Luxembourg norm which is equivalent to the Orlicz norm
(see \cite{bib:[K.M-W]} p 227 Th 7.5.4). It is convenient for us to work
with this norm (instead of the Orlicz norm). The space $L^{e}=\{f:\left\Vert
f\right\Vert _{(e)}<\infty \}$ is the Orlicz space with respect to the
Lebesgue measure. Notice that if we take $e_{p}(x)=\left\vert x\right\vert
^{p},p>1,$ then $\left\Vert f\right\Vert _{(e_{p})}$ is the usual $L^{p}$
norm and the corresponding Orlicz space is the standard $L^{p}$ space.
Another example is given at the end of this section.

\begin{remark}
\label{U}Let $u_{l}(x)=(1+\left\vert x\right\vert )^{-l}.$ As a consequence
of (\ref{Y}) ii), for every $l>d$ one has
\begin{equation*}
\left\Vert u_{l}\right\Vert _{e}\leq (e(1)\left\Vert u_{l}\right\Vert
_{1})\vee 1<\infty .
\end{equation*}%
Indeed (\ref{Y}) ii) implies that for $t\leq 1$ one has $e(t)\leq e(1)t.$
For $c>(e(1)\left\Vert u_{l}\right\Vert _{1})\vee 1$ one has $\frac{1}{c}%
u_{l}(x)\leq u_{l}(x)\leq 1$ so that
\begin{equation*}
\int e(\frac{1}{c}u_{l}(x))dx\leq \frac{e(1)}{c}\int u_{l}(x)dx=\frac{e(1)}{c%
}\left\Vert u_{l}\right\Vert _{1}\leq 1.
\end{equation*}
\end{remark}

One defines the conjugate of $e$ by $e_{\ast }(s)=\inf \{st-e(t):t\in \R\}$
and this is also a Young function so the corresponding Luxembourg norm $%
\left\Vert f\right\Vert _{(e_{\ast })}$ is given by (\ref{O1}) with $e$
replaced by $e_{\ast }$. Then we have the following H\"{o}lder inequality%
\begin{equation}
\left\vert \int fg(x)dx\right\vert \leq 2\left\Vert f\right\Vert
_{(e)}\left\Vert g\right\Vert _{(e_{\ast })}.  \label{O2}
\end{equation}%
See \cite{bib:[K.M-W]} p 215 Th 7.2.1 (the factor $2$ does not appear in Th
7.2.1 but there in the right hand side of the inequality one has the Orlicz
norm of $g.$ Using the equivalence between the Orlicz and the Luxembourg
norm we replace the Orlicz norm by $2\left\Vert g\right\Vert _{(e_{\ast
})}). $

If $e$ satisfies the $\Delta _{2}$ condition (that is (\ref{Y}) i) above)
then $L^{e}$ is reflexive (see \cite{bib:[K.M-W]} p 234 Th 7.7.1). In
particular, in this case, any bounded subset of $L^{e}$ is weakly relative
compact.

We also define $e^{-1}(a)=\sup \{c:e(c)\leq a\}$ and
\begin{equation}
\phi _{e}(r)=\frac{1}{e^{-1}(\frac{1}{r})}\quad and\quad \beta _{e}(R)=\frac{%
R}{e^{-1}(R)}=R\phi _{e}(\frac{1}{R}).  \label{O3}
\end{equation}

\begin{remark}
\label{Increasing} The function $\phi _{e}$ is the ``fundamental function''
of $L^{e}$ equipped with the Luxembourg norm (see \cite{bib:[BS]} Lemma 8.17
pg 276). In particular $\frac{1}{r}\phi _{e}(r)$ is decreasing (see \cite%
{bib:[BS]} Corollary 5.2 pg 67). It follows that $\beta _{e}$ is increasing.
For the sake of completeness we give here the argument. Indeed, if $a>1$
then $e(ax)\geq ae(x)$ so that $ax\geq e^{-1}(ae(x)).$ Taking $y=e(x)$ we
obtain $ae^{-1}(y)\geq e^{-1}(ay)$ which gives
\begin{equation*}
\beta _{e}(ay)=\frac{ay}{e^{-1}(ay)}\geq \frac{ay}{ae^{-1}(y)}=\beta _{e}(y).
\end{equation*}
\end{remark}

We introduce now the norms%
\begin{equation}
\left\Vert f\right\Vert _{k,(e)}=\sum_{0\leq \left\vert \alpha \right\vert
\leq k}\left\Vert \partial _{\alpha }f\right\Vert _{(e)}\quad and\quad
\left\Vert f\right\Vert _{k,\infty }=\sum_{0\leq \left\vert \alpha
\right\vert \leq k}\left\Vert \partial _{\alpha }f\right\Vert _{\infty }.
\label{O4}
\end{equation}

We denote%
\begin{equation*}
W^{k,e}=\{f:\left\Vert f\right\Vert _{k,(e)}<\infty \}\qquad and\qquad
W^{k,\infty }=\{f:\left\Vert f\right\Vert _{k,\infty }<\infty \}.
\end{equation*}

For a multi index $\gamma $ we denote $x^{\gamma
}=\prod_{i=1}^{d}x_{i}^{\gamma _{i}}$ and for two multi indexes $\alpha
,\gamma $ we denote $f_{\gamma ,\alpha }$ the function $f_{\alpha ,\gamma
}(x)=x^{\gamma }\partial _{\alpha }f(x).$ Then we consider the norm%
\begin{equation}
\left\Vert f\right\Vert _{k,l,(e)}=\sum_{0\leq \left\vert \gamma \right\vert
\leq l}\sum_{0\leq \left\vert \alpha \right\vert \leq k}\left\Vert f_{\gamma
,\alpha }\right\Vert _{(e)}\qquad and\qquad W^{k,l,e}=\{f:\left\Vert
f\right\Vert _{k,l,(e)}<\infty \}.  \label{O4a}
\end{equation}

We stress that in $\|\cdot \|_{k,l,(e)}$ the first index $k$ is related to
the order of the derivatives which are involved while the second index $l$
is connected to the power of the polynomial multiplying the function and its
derivatives up to order $k$.

\smallskip

We consider the following distances between two measures $\mu ,\nu \in
\mathcal{M}.$ For $k\in \N$%
\begin{equation}
d_{k}(\mu ,\nu )=\sup \{\left\vert \int \phi d\mu -\int \phi d\nu
\right\vert :\phi \in C^{\infty }(\R^{d}),\left\Vert \phi \right\Vert
_{k,\infty }\leq 1\}.  \label{O6}
\end{equation}%
Notice that $d_{0}$ is the total variation distance and $d_{1}$ is the
bounded variation distance. The Wasserstein distance (which is more popular)
is $d_{W}(\mu ,\nu )=\sup \{\left\vert \int \phi d\mu -\int \phi d\nu
\right\vert :\phi \in C^{1}(\R^{d}),\left\Vert \nabla \phi \right\Vert
_{\infty }\leq 1\}$ so $d_{1}(\mu ,\nu )\leq d_{W}(\mu ,\nu ).$ It follows
that all the results proved with respect to $d_{1}$ will be a fortiori true
for $d_{W}.$ The distances $d_{k}$ with $k\geq 2$ are less often used. We
mention however that people working in approximation theory (for diffusion
process for example - see \cite{bib:[TT]} or \cite{bib:[VN]}) use such
distances in an implicit way: indeed, they study the speed of convergence of
certain schemes but they are able to obtain their estimates for test
functions $f\in C^{k}$ with $k$ sufficiently large - so $d_{k}$ comes on.

Let $q,k\in \N$ and $m\in \N_{\ast }.$ For $\mu \in \mathcal{M}$ and for a
sequence $\mu _{n}\in \mathcal{M}_{a},n\in \N$ we define
\begin{equation}
\pi _{q,k,m,e}(\mu ,(\mu _{n})_{n})=\sum_{n=0}^{\infty }2^{n(q+k)}\beta
_{e}(2^{nd})d_{k}(\mu ,\mu _{n})+\sum_{n=0}^{\infty }\frac{1}{2^{2nm}}%
\left\Vert p_{\mu _{n}}\right\Vert _{2m+q,2m,(e)}.  \label{O7}
\end{equation}%
Moreover we define
\begin{equation}
\rho _{q,k,m,e}(\mu )=\inf \pi _{q,k,m,e}(\mu ,(\mu _{n})_{n})  \label{011}
\end{equation}%
with the infimum taken over all the sequences $\mu _{n}\in \mathcal{M}%
_{a},n\in \N.$ We define%
\begin{equation}
\mathcal{S}_{q,k,m,e}=\{\mu \in \mathcal{M}:\rho _{q,k,m,e}(\mu )<\infty \}.
\label{O13}
\end{equation}%
It is easy to check that $\rho _{q,k,m,e}$ is a norm on $\mathcal{S}%
_{q,k,m,e}.$

The main result in this section is the following.

\begin{theorem}
\label{2C} Let $q,k\in \N,m\in \N_{\ast }$ and let $e\in \mathcal{E}$.

i) Take $q=0.$ Then%
\begin{equation*}
\mathcal{S}_{0,k,m,e}\subset L^{e}
\end{equation*}%
in the sense that if $\mu \in \mathcal{S}_{0,k,m,e}$ then $\mu $ is
absolutely continuous and the density $p_{\mu }$ belongs to $L^{e}.$
Moreover there exists a universal constant $C$ such that
\begin{equation*}
\left\Vert p_{\mu }\right\Vert _{L^{e}}\leq C\rho _{0,k,m,e}(\mu ).
\end{equation*}

ii) Take $q\geq 1.$ Then
\begin{equation*}
\mathcal{S}_{q,k,m,e}\subset W^{q,e}\qquad and\qquad \left\Vert p_{\mu
}\right\Vert _{W^{q,e}}\leq C\rho _{q,k,m,e}(\mu ),\quad \mu \in \mathcal{S}%
_{q,k,m,e}.
\end{equation*}
\end{theorem}

The proof of this theorem is given at the end of Section \ref{sect-Hermite}.

It may be cumbersome to check that $\mu \in \mathcal{S}_{0,k,m,e}$ so we
give a sufficient condition which seems to be more clear and easier to
verify. We define
\begin{equation*}
M_{m,q,e}(R)=\{\mu \in \mathcal{M}_{a}\ :\ \left\Vert p_{\mu }\right\Vert
_{2m+q,2m,(e)}\leq R\}.
\end{equation*}%
For $a>1$ we denote%
\begin{equation*}
L_{a}(R)=R(\ln R)^{a}
\end{equation*}%
and we consider the following hypothesis.

\smallskip

\textbf{Hypothesis $H_{q}(k,m,e)$. } \emph{For $q,k\in \N$ , $m\in \N_{\ast }$
and} $e\in \mathcal{E}$ \emph{there exists $a>1$ such that%
\begin{equation*}
\overline{\lim }_{R\rightarrow \infty }\frac{L_{a}(R)^{1+\frac{k+q}{2m}%
}\beta _{e}(L_{a}(R)^{\frac{d}{2m}})}{R}d_{k}(\mu ,M_{m,q,e}(R))<\infty .
\end{equation*}%
}

\medskip

We define
\begin{equation*}
B_{q}(k,m,e)=\{\mu \in \mathcal{M}:\text{ }H_{q}(k,m,e)\text{ holds for }\mu
\}.
\end{equation*}

For $\alpha ,\gamma \geq 0$ we denote by $\mathcal{E}_{\alpha ,\gamma }$ the
class of the Young functions $e\in \mathcal{E}$ such that
\begin{equation}
0<\underline{\lim }_{R\rightarrow \infty }\frac{\beta _{e}(R)}{R^{\alpha
}(\ln R)^{\gamma }}\leq \overline{\lim }_{R\rightarrow \infty }\frac{\beta
_{e}(R)}{R^{\alpha }(\ln R)^{\gamma }}<\infty .  \label{RR}
\end{equation}%
The examples we have in mind fit in this class.

Our criterion is the following.

\begin{theorem}
\label{2E} i) Let $q,k\in \N$, $m\in \N_{\ast }$ and $e\in \mathcal{E}$. Then%
\begin{equation*}
B_{q}(k,m,e)\subset \mathcal{S}_{q,k,m,e}\subset W^{q,e}.
\end{equation*}

ii) Suppose that $e\in \mathcal{E}_{\alpha ,\gamma }$ with $0\leq \alpha <%
\frac{2m+k+q}{d(2m-1)}$ and $\gamma \geq 0.$ Then%
\begin{equation*}
W^{q+1,2m,e}\subset B_{q}(k,m,e)\subset \mathcal{S}_{q,k,m,e}\subset W^{q,e},
\end{equation*}
the first inclusion holding for $m>d/2$.
\end{theorem}

\textbf{Proof.} i) The proof of $B_{q}(k,m,e)\subset \mathcal{S}_{q,k,m,e}$
is given in Lemma \ref{BALANCE} and $\mathcal{S}_{q,k,m,e}\subset W^{q,e}$
is Theorem \ref{2C}.

ii) The inclusion $W^{q+1,2m,e}\subset B_{q}(k,m,e)$ is proved in Corollary %
\ref{Reciproc}. $\square $

\begin{remark}
The above criterion involves a lot of parameters and it is not easy to
understand at a first glance which is the significance of each of them. So
we try to give a first interpretation on their meaning. Our aim is to prove
that the measure $\mu $ has a certain regularity and we want to do this by
approximating it by some regular measures $\nu \in M_{m,q,e}(R).$ The first
parameter which we fix is $q.$ It represents the order of regularity that we
hope to obtain. If $q=0$ this means that we want to prove that $\mu $ is
just absolutely continuous with respect to the Lebesgue measure and the
density is in some Orlicz space defined by $e.$ For example if $%
e(t)=\left\vert t\right\vert ^{p}$ this is the $L^{p}$ space. If $q\geq 1$
then our aim is to obtain more regularity, namely to prove that the density
is in the Sobolev space of order $q.$ The second parameter is $k.$ It
characterizes the distance in which we estimate the approximation error.
Once $q,e,k$ are chosen it remains $m.$ The choice of $m$ is different. For
the other parameters the choice comes from our decision to treat a problem
or another: they are involved in the definition of the problem we treat. On
the contrary, $m$ is a free parameter which we choose for technical reasons:
for example in some concrete situations it is suitable to chose $m$ very
large - such that $\frac{1}{m}\leq \varepsilon $ for some $\varepsilon >0$
for example. In the following it will become clear that $m$ represents the
number of integration by parts that we use.
\end{remark}

The hypothesis $H_{q}(k,m,e)$ seems difficult to check because of the
function $\beta _{e}$ which is involved there in. So we give two significant
examples.

\medskip

\textbf{Example 1.} 
\emph{Let $p>1$ and $e_{p}(t)=\left\vert t\right\vert ^{p}.$ The
corresponding Orlicz space is the standard $L^{p}$ space. And we have $\beta
_{e_{p}}(t)=\left\vert t\right\vert ^{1/p_{\ast }}$ with $p_{\ast }$ the
conjugate of $p.$ So $e_{p}\in \mathcal{E}_{1/p_{\ast },0}.$ Our hypothesis
is:
\begin{equation*}
H_{q}(k,m,e_{p}):\quad \exists a>1\text{ s.t. }\overline{\lim }%
_{R\rightarrow \infty }R^{\frac{q+k+d/p_{\ast }}{2m}}(\ln R)^{a(1+\frac{%
q+k+d/p_{\ast }}{2m})}d_{k}(\mu ,M_{m,q,e}(R))<\infty .
\end{equation*}%
In this hypothesis the dimension $d$ is still present. But its contribution
is very small when $p$ is close to one. Then Theorem \ref{2E} reads as
follows. If $\mu $ satisfies $H_{q}(k,m,e_{p})$ then $\mu (dx)=f(x)dx$ with $%
f\in W^{q,2m,p}.$ Moreover, if $p_{\ast }>\frac{d(2m-1)}{2m+k+q}$ and $f\in
W^{q+1,p}$ then the measure $\mu (dx)=f(x)dx$ satisfies $H_{q}(k,m,e_{p}).$
}

\medskip

\textbf{Example 2.} 
\emph{Set $e_{\log }(t)=(1+\left\vert t\right\vert )\ln (1+\left\vert
t\right\vert ).$ It is easy to check that $e_{\log }$ is a Young function
which satisfies the property $\Delta _{2}.$ The corresponding Orlicz space
is the so called LlogL space of Zigmund, i.e the spaces of the
functions $f$ such that $\int \left\vert f(x)\right\vert \ln
^{+}\left\vert f(x)\right\vert dx<\infty$ (see \cite{bib:[BS]}).
Let us check that $e_{\log }\in \mathcal{E}_{0,1}.$ We denote $%
e_{a}(t)=at\ln (at)$ and we notice that for $t$ large we have $e_{1}(t)\leq
e_{\log }(t)\leq e_{2}(t)$ (in particular $L^{e_{\log }}$ is the space of
the functions which have finite entropy) Then $e_{2}^{-1}(t)\leq e_{\log
}^{-1}(t)\leq e_{1}^{-1}(t)$ so that
\begin{equation*}
\frac{t}{e_{1}^{-1}(t)}\leq \beta _{e_{\log }}(t)\leq \frac{t}{e_{2}^{-1}(t)}.
\end{equation*}%
Using the change of variable $R=e_{a}(t)$ one obtains%
\begin{equation*}
\lim_{R\rightarrow \infty }\frac{R}{e_{1}^{-1}(R)\ln R}=\lim_{t\rightarrow
\infty }\frac{e_{a}(t)}{t\ln e_{a}(t)}=a.
\end{equation*}%
This proves that
\begin{equation*}
1\leq \underline{\lim }_{t\rightarrow \infty }\frac{\beta _{e_{\log }}(t)}{%
\ln t}\leq \overline{\lim }_{t\rightarrow \infty }\frac{\beta _{e_{\log }}(t)%
}{\ln t}\leq 2.
\end{equation*}%
Our hypothesis is:
\begin{equation}
H_{q}(k,m,e_{\log }):\quad \exists a>1\text{ s.t. }\overline{\lim }%
_{R\rightarrow \infty }R^{\frac{q+k}{2m}}(\ln R)^{a(1+\frac{q+k}{2m}%
)+1}d_{k}(\mu ,M_{m,q,e_{\log }}(R))<\infty .  \label{O11}
\end{equation}%
So the dimension $d$ does no more appear. Then Theorem \ref{2E} reads as
follows. If $\mu $ satisfies $H_{q}(k,m,e_{\log })$ then $\mu (dx)=p(x)dx$
with $p\in W^{q,e_{\log }}.$ And if $p\in W^{q+1,2m,e_{\log }}$ then the
measure $\mu (dx)=p(x)dx$ satisfies $H_{q}(k,m,e_{\log }).$
}

\begin{remark}
One may consider a further step and take $e_{\log \log }(t)=(1+\left\vert
t\right\vert )\ln (1+\ln (1+\left\vert t\right\vert )).$ In this case $\beta
_{e_{\log \log }}(t)\leq \ln \ln t$ but there will be no significant
improvement: we obtain $(\ln R)^{a(1+\frac{q+k}{2m})}\ln \ln R$ instead of $%
(\ln R)^{a(1+\frac{q+k}{2m})+1}.$
\end{remark}

\begin{remark}
Having in mind these examples we conclude that one may ask for two different
questions. a) One may just want to prove that $\mu $ is absolutely
continuous with respect to the Lebesgue measure. b) One wants to obtain
estimates of the density in some given norm (associated to some Young
function $e$ - for example in $L^{p}).$ If the question is just a) then one
has to go directly to Example 2 and to use the Young function $e_{\log }.$
Because the hypothesis (\ref{O11}) is the minimal one (in our approach at
least).
\end{remark}


\subsection{Hermite expansions and density estimates}

\label{sect-Hermite}

We begin with a re-view of some basic properties of Hermite polynomials and
functions. The Hermite polynomials on $\R$ are defined by%
\begin{equation*}
H_{n}(t)=(-1)^{n}e^{t^{2}}\frac{d^{n}}{dt}e^{-t^{2}},\quad n=0,1,...
\end{equation*}%
They are orthogonal with respect to $e^{-t^{2}}dt.$ We denote the $L^{2}$
normalized Hermite functions by
\begin{equation*}
h_{n}(t)=(2^{n}n!\sqrt{\pi })^{-1/2}H_{n}(t)e^{-t^{2}/2}
\end{equation*}%
and we have%
\begin{equation*}
\int_{R}h_{n}(t)h_{m}(t)dt=(2^{n}n!\sqrt{\pi })^{-1}%
\int_{R}H_{n}(t)H_{m}(t)e^{-t^{2}}dt=\delta _{n,m}.
\end{equation*}%
The Hermite functions form an orthonormal basis in $L^{2}(R).$ For a multi
index $\alpha =(\alpha _{1},...,\alpha _{d})\in \N^{d}$ we define the $d$%
-dimensional Hermite function
\begin{equation*}
\mathcal{H}_{\alpha }(x):=\prod_{i=1}^{d}h_{\alpha _{i}}(x_{i}),\quad
x=(x_{1},...,x_{d}).
\end{equation*}%
The $d$-dimensional Hermite functions form an orthonormal basis in $%
L^{2}(\R^{d}).$ This corresponds to the chaos decomposition in dimension $d$
(but the notation we gave above is slightly different from the one used in
probability; see \cite{bib:[N]}, \cite{bib:[Sa]} and \cite{bib:[M]}, where
Hermite polynomials are used. One may come back by a renormalization). The
Hermite functions are the eigenvectors of the Hermite operator $D=-\Delta
+\left\vert x\right\vert ^{2}$ (with $\Delta $ the Laplace operator) and one
has%
\begin{equation}
D\mathcal{H}_{\alpha }=(2\left\vert \alpha \right\vert +d)\mathcal{H}%
_{\alpha }\quad with\quad \left\vert \alpha \right\vert =\alpha
_{1}+...+\alpha _{d}.  \label{Her1}
\end{equation}%
We denote $W_{n}=span\{\mathcal{H}_{\alpha }:\left\vert \alpha \right\vert
=n\}$ and we have $L^{2}(\R^{d})=\oplus _{n=0}^{\infty }W_{n}$.

For a function $\Phi :\R^{d}\times \R^{d}\rightarrow \R$ and a function $%
f:\R^{d}\rightarrow \R$ we use the notation%
\begin{equation*}
\Phi \ast f(x)=\int_{\R^{d}}\Phi (x,y)f(y)dy.
\end{equation*}%
We denote by $J_{n}$ the orthogonal projection on $W_{n}$ and we have
\begin{equation}
J_{n}v(x)=\mathcal{H}_{n}\ast v(x)\quad with\quad \mathcal{H}%
_{n}(x,y):=\sum_{\left\vert \alpha \right\vert =n}\mathcal{H}_{\alpha }(x)%
\mathcal{H}_{\alpha }(y).  \label{Her2}
\end{equation}%
Moreover, we consider a function $a:\R_{+}\rightarrow \R$ whose support is
included in $[\frac{1}{4},4]$\ and we define%
\begin{equation*}
\mathcal{H}_{n}^{a}(x,y)=\sum_{j=0}^{\infty }a(\frac{j}{4^{n}})\mathcal{H}%
_{j}(x,y) =\sum_{j=4^{n-1}+1}^{4^{n+1}-1}a(\frac{j}{4^{n}})\mathcal{H}%
_{j}(x,y),\quad x,y\in \R^d,
\end{equation*}%
the last equality being a consequence of the support property of the
function $a.$

The following estimate is a crucial point in our approach. It has been
proved in \cite{bib:[E]}, \cite{bib:[Dz]} and then in \cite{bib:[PY]}. We
refer to Corollary 2.3, inequality (2.17), in \cite{bib:[PY]}.

\begin{theorem}
\label{5}Let $a:\R_{+}\rightarrow \R_{+}$ be a non negative $C^{\infty }$
function with the support included in $[\frac{1}{4},4].$ We denote $%
\left\Vert a\right\Vert _{l}=\sum_{i=0}^{l}\sup_{t\geq 0}\left\vert
a^{(i)}(t)\right\vert .$ For every multi-index $\alpha $ and every $k\in \N$
there exists a constant $C_{k}$ (depending on $k,\alpha ,d)$ such that for
every $n\in \N$ and every $x,y\in \R^{d}$%
\begin{equation}
\left\vert \frac{\partial ^{\left\vert \alpha \right\vert }}{\partial
x^{\alpha }}\mathcal{H}_{n}^{a}(x,y)\right\vert \leq C_{k}\left\Vert
a\right\Vert _{k}\frac{2^{n(\left\vert \alpha \right\vert +d)}}{%
(1+2^{n}\left\vert x-y\right\vert )^{k}}.  \label{Her3}
\end{equation}
\end{theorem}

Following the ideas in \cite{bib:[PY]} we consider a function $%
a:\R_{+}\rightarrow \R_{+} $ of class $C_{b}^{\infty }$ with the support
included in $[\frac{1}{4},4]$ and such that $a(t)+a(4t)=1$ for $t\in \lbrack
\frac{1}{4},1].$ We may construct $a$ in the following way: we take a
function $a:[0,1]\rightarrow \R_{+}$ with $a(t)=0$ for $t\leq \frac{1}{4}$
and $a(1)=1.$ We may choose $a$ such that $a^{(l)}(\frac{1}{4}%
)=a^{(l)}(1-)=0 $ for every $l\in \N$. Then we define $a(t)=1-a(\frac{t}{4})$
for $t\in \lbrack 1,4]$ and $a(t)=0$ for $t\geq 4.$ This is the function we
will use in the following$.$ Notice that $\ a$ has the property:%
\begin{equation}
\sum_{n=0}^{\infty }a(\frac{t}{4^{n}})=1\quad \forall t\geq 1.  \label{Her0}
\end{equation}%
In order to check the above equality we fix $n_{t}$ such that $%
4^{n_{t}-1}\leq t<4^{n_{t}}$ and we notice that $a(\frac{t}{4^{n}})=0$ if $%
n\notin \{n_{t}-1,n_{t}\}.$ So $\sum_{n=0}^{\infty }a(\frac{t}{4^{n}}%
)=a(4s)+a(s)=1$ with $s=t/4^{n_{t}}\in \lbrack \frac{1}{4},1].$ In the
following we fix a function $a$ and the constants in our estimates will
depend on $\left\Vert a\right\Vert _{l}$ for some fixed $l.$ Using this
function we obtain the following representation formula:

\begin{proposition}
\label{6}For every $f\in L^{2}(\R^{d})$
\begin{equation*}
f=\sum_{n=0}^{\infty }\mathcal{H}_{n}^{a}\ast f
\end{equation*}%
the series being convergent in $L^{2}(\R^{d}).$
\end{proposition}

\textbf{Proof}. We fix $N$ and we denote
\begin{equation*}
S_{N}^{a}=\sum_{n=1}^{N}\mathcal{H}_{n}^{a}\ast f,\quad
S_{N}=\sum_{j=1}^{4^{N}}\mathcal{H}_{j}\ast f\quad and\quad
R_{N}^{a}=\sum_{j=4^{N}+1}^{4^{N+1}}(\mathcal{H}_{j}\ast f)a(\frac{j}{4^{N+1}%
}).
\end{equation*}%
Let $j\leq 4^{N+1}.$ For $n\geq N+2$ one has $a(\frac{j}{4^{n}})=0.$ So
using (\ref{Her0}) we obtain $\sum_{n=1}^{N}a(\frac{j}{4^{n}}%
)=\sum_{n=1}^{\infty }a(\frac{j}{4^{n}})-a(\frac{j}{4^{N+1}})=1-a(\frac{j}{%
4^{N+1}}).$ And for $j\leq 4^{N}$ one has $a(\frac{j}{4^{N+1}})=0.$\ It
follows that%
\begin{eqnarray*}
S_{N}^{a} &=&\sum_{n=1}^{N}\sum_{j=0}^{\infty }a(\frac{j}{4^{n}})\mathcal{H}%
_{j}\ast f=\sum_{n=1}^{N}\sum_{j=0}^{4^{N+1}}a(\frac{j}{4^{n}})\mathcal{H}%
_{j}\ast f=\sum_{j=0}^{4^{N+1}}(\mathcal{H}_{j}\ast f)\sum_{n=1}^{N}a(\frac{j%
}{4^{n}}) \\
&=&\sum_{j=0}^{4^{N+1}}\mathcal{H}_{j}\ast f-\sum_{j=4^{N}+1}^{4^{N+1}}(%
\mathcal{H}_{j}\ast f)a(\frac{j}{4^{N+1}})=S_{N+1}-R_{N}^{a}.
\end{eqnarray*}%
One has $S_{N}\rightarrow f$ in $L^{2}$ and $\left\Vert R_{N}^{a}\right\Vert
_{2}\leq \left\Vert a\right\Vert _{\infty
}\sum_{j=4^{N}+1}^{4^{N+1}}\left\Vert \mathcal{H}_{j}\ast f\right\Vert
_{2}\rightarrow 0$ so the proof is completed. $\square $

\smallskip

We will need the following lemma concerning properties of the Luxembourg
norms.

\begin{lemma}
Let $\rho \geq 0$ be a measurable function. Then for every measurable
function $f$
\begin{equation}
\left\Vert \rho \ast f\right\Vert _{(e)}\leq \left\Vert \rho \right\Vert
_{1}\left\Vert f\right\Vert _{(e)}.  \label{Oo2}
\end{equation}
\end{lemma}

\textbf{Proof}. Let $c=m\left\Vert f\right\Vert _{(e)}$ with $m=\left\Vert
\rho \right\Vert _{1}=\int \rho (x-y)dy.$ Since $e$ is convex we obtain%
\begin{eqnarray*}
\int e(\frac{1}{c}(\rho \ast f)(x))dx &=&\int e(\int \frac{\rho (x-y)}{m}%
\times \frac{m}{c}f(y)dy)dx \\
&\leq &\int dx\int \frac{\rho (x-y)}{m}\times e(\frac{m}{c}f(y))dy \\
&=&\int e(\frac{m}{c}f(y))\int \frac{\rho (x-y)}{m}dxdy=\int e(\frac{m}{c}%
f(y))dy \\
&=&\int e(\frac{1}{\left\Vert f\right\Vert _{(e)}}f(y))dy\leq 1
\end{eqnarray*}%
and this means that $\left\Vert \rho \ast f\right\Vert _{(e)}\leq
c=\left\Vert \rho \right\Vert _{1}\left\Vert f\right\Vert _{(e)}.$ $\square $

\begin{lemma}
Let $\rho _{n,p}(z)=(1+2^{n}\left\vert z\right\vert )^{-p}$ with $p>d$ and $%
e\in \mathcal{E}.$ There exists a constant $C_{p}$ depending on $p$ and $d$
such that
\begin{equation}
\left\Vert \rho _{n,p}\right\Vert _{(e)}\leq \frac{1}{e^{-1}(\frac{1}{C_{p}}%
2^{nd})}  \label{Oo1}
\end{equation}%
In particular, for $p=d+1$ there exists a constant $C$ depending on $d$ and
on the doubling constant of $e$ such that (with $\phi _{e}$ defined in (\ref%
{O3}))%
\begin{equation}
\left\Vert \rho _{n,d+1}\right\Vert _{(e)}\leq \frac{C}{e^{-1}(2^{nd})}%
=C2^{-nd}\beta _{e}(2^{nd})=C\phi _{e}(\frac{1}{2^{nd}})  \label{Oo1a}
\end{equation}
\end{lemma}

\textbf{Proof}. Let $c>0.$ We pass in polar coordinates first and we use the
change of variable $s=2^{n}r$ then and we obtain%
\begin{eqnarray*}
\int_{\R^{d}}e(\frac{1}{c}\rho _{n,p}(z))dz &=&A_{d}\int_{0}^{\infty
}r^{d-1}e(\frac{1}{c}\times \frac{1}{(1+2^{n}r)^{p}})dr \\
&=&2^{-nd}A_{d}\int_{0}^{\infty }s^{d-1}e(\frac{1}{c}\times \frac{1}{%
(1+s)^{p}})ds
\end{eqnarray*}%
where $A_{d}$ is the surface of the unit sphere in $\R^{d}.$ Using the
property (\ref{Y}) ii) we upper bound the above term by%
\begin{equation*}
2^{-nd}e(\frac{1}{c})A_{d}\int_{0}^{\infty }s^{d-1}\times \frac{1}{(1+s)^{p}}%
ds=C_{p}2^{-nd}e(\frac{1}{c}).
\end{equation*}%
In order to prove that $\left\Vert \rho _{n,p}\right\Vert _{(e)}\leq c$ we
have to check that $\int_{\R^{d}}e(\frac{1}{c}\rho _{n,p}(z))dz\leq 1.$ In
view of the above inequalities it suffices that $e(\frac{1}{c})\leq
2^{nd}/C_{p}$ that is $c\geq 1/e^{-1}(2^{nd}/C_{p}).$ $\square $

\begin{proposition}
\label{7a} Let $\alpha $ be a multi index.

i) There exists a universal constant $C$ (depending on $\alpha ,d$ and $e)$
such that%
\begin{eqnarray}
a)\quad \left\Vert \partial _{\alpha }\mathcal{H}_{n}^{a}\ast f\right\Vert
_{(e)} &\leq &C\left\Vert a\right\Vert _{d+1}\times 2^{n\left\vert \alpha
\right\vert }\left\Vert f\right\Vert _{(e)},  \label{Oo3} \\
b)\quad \left\Vert \partial _{\alpha }\mathcal{H}_{n}^{a}\ast f\right\Vert
_{\infty } &\leq &C\left\Vert a\right\Vert _{d+1}\times 2^{n\left\vert
\alpha \right\vert }\beta _{e}(2^{nd})\left\Vert f\right\Vert _{(e^{\ast })}
\notag
\end{eqnarray}

ii) Let $m\in \N_{\ast }.$ There exists a universal constant $C$ (depending
on $\alpha ,m,d$ and $e)$ such that
\begin{equation}
\left\Vert \mathcal{H}_{n}^{a}\ast \partial _{\alpha }f\right\Vert
_{(e)}\leq \frac{C\left\Vert a\right\Vert _{d+1}^{2}}{4^{nm}}\left\Vert
f\right\Vert _{2m+\left\vert \alpha \right\vert ,2m,(e)}  \label{Oo4}
\end{equation}

iii) Let $k\in \N.$ There exists a universal constant $C$ (depending on $%
\alpha ,k,d$ and $e)$ such that
\begin{equation}
\left\Vert \mathcal{H}_{n}^{a}\ast \partial _{\alpha }(f-g)\right\Vert
_{(e)}\leq C\left\Vert a\right\Vert _{d+1}\times 2^{n(\left\vert \alpha
\right\vert +k)}\beta (2^{nd})d_{k}(\mu _{f},\mu _{g})  \label{Oo5}
\end{equation}%
iv) For $f,g\in W^{k,e_{\ast }}$ we define
\begin{equation*}
d_{k,e_{\ast }}(\mu _{f},\mu _{g})=\sup \{\left\vert \int \phi d\mu
_{f}-\int \phi d\mu _{g}\right\vert :\phi \in W^{k,e_{\ast }},\left\Vert
\phi \right\Vert _{(e_{\ast })}\leq 1\}.
\end{equation*}%
Let $k\in \N.$ There exists a universal constant $C$ (depending on $\alpha
,k,d$ and $e)$ such that%
\begin{equation}
\left\Vert \mathcal{H}_{n}^{a}\ast \partial _{\alpha }(f-g)\right\Vert
_{(e)}\leq C\left\Vert a\right\Vert _{d+1}\times 2^{n(\left\vert \alpha
\right\vert +k)}d_{k,e_{\ast }}(\mu _{f},\mu _{g}).  \label{Oo6}
\end{equation}
\end{proposition}

\begin{remark}
Let us try to explain the gain with respect to the standard projection on
the basis of Hermite functions. We have to compare $\mathcal{H}_{n}^{a}$
with $\mathcal{H}_{4^{n}}$ because $\mathcal{H}_{n}^{a}$ is a mixture of
projections between $4^{n-1}$ and $4^{n+1}.$ Suppose that we work in $L^{p}$
so $e_{p}(x)=\left\vert x\right\vert ^{p}.$ In this case $\beta_e
(t)=C\times t^{1/p_{\ast }}$ (with $p_{\ast }$ the conjugate of $p)$ so the
inequality (\ref{Oo5}), with $k=1,$ reads%
\begin{equation}
\left\Vert \mathcal{H}_{n}^{a}\ast \partial _{\alpha }(f-g)\right\Vert
_{p}\leq C2^{n(\left\vert \alpha \right\vert +1+\frac{d}{p_{\ast }}%
)}d_{1}(\mu _{f},\mu _{g})  \label{Comp1}
\end{equation}%
In some notes (which are not reported here) we were able to prove that
\begin{equation}
\left\Vert \mathcal{H}_{4^{n}}\ast \partial _{\alpha }(f-g)\right\Vert
_{2}\leq C2^{n(\left\vert \alpha \right\vert +1+\frac{d}{2})}d_{1}(\mu
_{f},\mu _{g}).  \label{Comp2}
\end{equation}%
If we take $p=2$ in (\ref{Comp1}) we obtain exactly (\ref{Comp2}) so there
is no gain. But we could not obtain the estimates in (\ref{Comp2}) for any $%
p,$ but only for $p=2.$ Here we are able to take any $p>1$ so $p_{\ast } $
may be taken arbitrary large - and this destroys the dimension $d.$
\end{remark}

\textbf{Proof}. i) Using (\ref{Her3})%
\begin{equation}
\left\vert \partial _{\alpha }\mathcal{H}_{n}^{a}\ast f(x)\right\vert \leq
C2^{n(\left\vert \alpha \right\vert +d)}\left\Vert a\right\Vert _{d+1}\int
\rho _{n,d+1}(x-y)\left\vert f(y)\right\vert dy.  \label{Oo5'}
\end{equation}%
Since $e$ is symmetric $e(\left\vert x\right\vert )=e(x)$ so that $%
\left\Vert f\right\Vert _{(e)}=\left\Vert \left\vert f\right\vert
\right\Vert _{(e)}.$ Moreover, if $0\leq f(x)\leq g(x)$ then $\left\Vert
f\right\Vert _{(e)}\leq \left\Vert g\right\Vert _{(e)}.$ Using these
properties and (\ref{Oo5'}) and (\ref{Oo2}) we obtain
\begin{eqnarray*}
\left\Vert \partial _{\alpha }\mathcal{H}_{n}^{a}\ast f\right\Vert _{(e)}
&=&\left\Vert \left\vert \partial _{\alpha }\mathcal{H}_{n}^{a}\ast
f\right\vert \right\Vert _{(e)}\leq C2^{n(\left\vert \alpha \right\vert
+d)}\left\Vert a\right\Vert _{d+1}\left\Vert \rho _{n,d+1}\ast \left\vert
f\right\vert \right\Vert _{(e)} \\
&\leq &C2^{n(\left\vert \alpha \right\vert +d)}\left\Vert a\right\Vert
_{d+1}\left\Vert \rho _{n,d+1}\right\Vert _{1}\left\Vert \left\vert
f\right\vert \right\Vert _{(e)}.
\end{eqnarray*}%
Using (\ref{Oo1a}) with $e(x)=\left\vert x\right\vert $ we obtain $%
\left\Vert \rho _{n,d+1}\right\Vert _{1}\leq C/2^{nd}.$ So we conclude that
\begin{equation*}
\left\Vert \partial _{\alpha }\mathcal{H}_{n}^{a}\ast f\right\Vert
_{(e)}\leq C\left\Vert a\right\Vert _{d+1}2^{n\left\vert \alpha \right\vert
}\left\Vert f\right\Vert _{(e)}
\end{equation*}%
so a) is proved. Again by (\ref{Oo5'})
\begin{align*}
\left\vert \partial _{\alpha }\mathcal{H}_{n}^{a}\ast f(x)\right\vert & \leq
C\left\Vert a\right\Vert _{d+1}2^{n(\left\vert \alpha \right\vert +d)}\int
\rho _{n,d+1}(x-y)\left\vert f(y)\right\vert dy \\
& \leq C\left\Vert a\right\Vert _{d+1}2^{n(\left\vert \alpha \right\vert
+d)}\left\Vert \rho _{n,d+1}\right\Vert _{(e)}\left\Vert f\right\Vert
_{(e_{\ast })},
\end{align*}
the second inequality being a consequence of the H\"{o}lder inequality.
Using (\ref{Oo1a}), b) is proved as well.

\smallskip

ii) We define the functions $a_{m}(t)=a(t)t^{-m}.$ Since $a(t)=0$ for $t\leq
\frac{1}{4}$ and for $t\geq 4$ we have $\left\Vert a_{m}\right\Vert
_{d+1}\leq C_{m,d}\left\Vert a\right\Vert _{d+1}.$ Moreover $D\mathcal{H}%
_{j}\ast v=(2j+d)\mathcal{H}_{j}\ast v$ so we obtain%
\begin{equation*}
\mathcal{H}_{j}\ast v=\frac{1}{2j}(D-d)\mathcal{H}_{j}\ast v.
\end{equation*}%
We denote $L_{m,\alpha }=(D-d)^{m}\partial _{\alpha }$ and we notice that $%
L_{m,\alpha }=\sum_{\left\vert \beta \right\vert \leq 2m}\sum_{\left\vert
\gamma \right\vert \leq 2m+\left\vert \alpha \right\vert }c_{\beta ,\gamma
}x^{\beta }\partial _{\gamma }$ where $c_{\beta ,\gamma }$ are universal
constants. It follows that there exists some universal constant $C$ such
that
\begin{equation}
\left\Vert L_{m,\alpha }f\right\Vert _{(e)}\leq C\left\Vert f\right\Vert
_{2m+\left\vert \alpha \right\vert ,2m,(e)}.  \label{Oo9}
\end{equation}%
We take now $v\in L^{e_{\ast }}$ and we write
\begin{eqnarray*}
\left\langle v,\mathcal{H}_{n}^{a}\ast (\partial _{\alpha }f)\right\rangle
&=&\left\langle \mathcal{H}_{n}^{a}\ast v,\partial _{\alpha }f\right\rangle
=\sum_{j=0}^{\infty }a(\frac{j}{4^{n}})\left\langle \mathcal{H}_{j}\ast
v,\partial _{\alpha }f\right\rangle \\
&=&\sum_{j=1}^{\infty }a(\frac{j}{4^{n}})\frac{1}{(2j)^{m}}\left\langle
(D-d)^{m}\mathcal{H}_{j}\ast v,\partial _{\alpha }f\right\rangle \\
&=&\frac{1}{2^{m}}\times \frac{1}{4^{nm}}\sum_{j=1}^{\infty }a_{m}(\frac{j}{%
4^{n}})\left\langle \mathcal{H}_{j}\ast v,L_{m,\alpha }f\right\rangle \\
&=&\frac{1}{2^{m}}\times \frac{1}{4^{nm}}\left\langle \mathcal{H}%
_{n}^{a_{m}}\ast v,L_{m,\alpha }f\right\rangle .
\end{eqnarray*}%
Using the decomposition form Proposition \ref{6}\ we write $L_{m,\alpha
}f=\sum_{j=0}^{\infty }\mathcal{H}_{j}^{a}\ast L_{m,\alpha }f.$ For $j<n-1$
and for $j>n+1$ we have $\left\langle \mathcal{H}_{n}^{a_{m}}\ast v,\mathcal{%
H}_{j}^{a}\ast L_{m,\alpha }f\right\rangle =0.$ So using H\"{o}lder's
inequality%
\begin{eqnarray*}
\left\vert \left\langle v,\mathcal{H}_{n}^{a}\ast (\partial _{\alpha
}f)\right\rangle \right\vert &\leq &\frac{1}{4^{nm}}\sum_{j=n-1}^{n+1}\left%
\vert \left\langle \mathcal{H}_{n}^{a_{m}}\ast v,\mathcal{H}_{j}^{a}\ast
L_{m,\alpha }f\right\rangle \right\vert \\
&\leq &\frac{1}{4^{nm}}\sum_{j=n-1}^{n+1}\left\Vert \mathcal{H}%
_{n}^{a_{m}}\ast v\right\Vert _{(e_{\ast })}\left\Vert \mathcal{H}%
_{j}^{a}\ast L_{m,\alpha }f\right\Vert _{(e)}.
\end{eqnarray*}%
Using the point i) a) with $\alpha $ the void index we obtain $\left\Vert
\mathcal{H}_{n}^{a_{m}}\ast v\right\Vert _{(e_{\ast })}\leq C\left\Vert
a_{m}\right\Vert _{d+1}\left\Vert v\right\Vert _{(e_{\ast })}\leq C\times
C_{m,d}\left\Vert a\right\Vert _{d+1}\left\Vert v\right\Vert _{(e_{\ast })}$%
. Moreover, we have $\left\Vert \mathcal{H}_{j}^{a}\ast L_{m,\alpha
}f\right\Vert _{(e)}\leq C\left\Vert a\right\Vert _{d+1}\left\Vert
L_{m,\alpha }f\right\Vert _{(e)}\leq C\left\Vert a\right\Vert
_{d+1}\left\Vert f\right\Vert _{2m+\left\vert \alpha \right\vert ,2m,(e)},$
the last inequality being a consequence of (\ref{Oo9}). We obtain%
\begin{equation*}
\left\vert \left\langle v,\mathcal{H}_{n}^{a}\ast (\partial _{\alpha
}f)\right\rangle \right\vert \leq \frac{C\left\Vert a\right\Vert _{d+1}^{2}}{%
4^{nm}}\left\Vert v\right\Vert _{e_{\ast }}\left\Vert f\right\Vert
_{2m+\left\vert \alpha \right\vert ,2m,(e)}
\end{equation*}%
and, since $L^{e}$ is reflexive, (\ref{Oo4}) is proved.

\smallskip

iii) We write
\begin{eqnarray*}
\left\vert \left\langle v,\mathcal{H}_{n}^{a}\ast (\partial _{\alpha
}(f-g))\right\rangle \right\vert &=&\left\vert \left\langle \mathcal{H}%
_{n}^{a}\ast v,\partial _{\alpha }(f-g)\right\rangle \right\vert =\left\vert
\left\langle \partial _{\alpha }\mathcal{H}_{n}^{a}\ast v,f-g)\right\rangle
\right\vert \\
&=&\left\vert \int \partial _{\alpha }\mathcal{H}_{n}^{a}\ast vd\mu
_{f}-\int \partial _{\alpha }\mathcal{H}_{n}^{a}\ast vd\mu _{g}\right\vert .
\end{eqnarray*}%
We use (\ref{Oo3}) b)\ and we obtain%
\begin{eqnarray*}
&&\left\vert \int \partial _{\alpha }\mathcal{H}_{n}^{a}\ast vd\mu _{f}-\int
\partial _{\alpha }\mathcal{H}_{n}^{a}\ast vd\mu _{g}\right\vert \leq
\left\Vert \partial _{\alpha }\mathcal{H}_{n}^{a}\ast v\right\Vert
_{k,\infty }d_{k}(\mu _{f},\mu _{g}) \\
&\leq &\left\Vert \mathcal{H}_{n}^{a}\ast v\right\Vert _{k+\left\vert \alpha
\right\vert ,\infty }d_{k}(\mu _{f},\mu _{g})\leq C\left\Vert a\right\Vert
_{d+1}2^{n(k+\left\vert \alpha \right\vert )}\beta _{e}(2^{nd})\left\Vert
v\right\Vert _{(e_{\ast })}d_{k}(\mu _{f},\mu _{g})
\end{eqnarray*}%
which implies (\ref{Oo5}).

\smallskip

iv) We use (\ref{Oo3}) a)\ and we obtain%
\begin{eqnarray*}
\left\vert \int \partial _{\alpha }\mathcal{H}_{n}^{a}\ast vd\mu _{f}-\int
\partial _{\alpha }\mathcal{H}_{n}^{a}\ast vd\mu _{g}\right\vert &\leq
&\left\Vert \partial _{\alpha }\mathcal{H}_{n}^{a}\ast v\right\Vert
_{k,(e_{\ast })}d_{k,e_{\ast }}(\mu _{f},\mu _{g}) \\
&\leq &\left\Vert \mathcal{H}_{n}^{a}\ast v\right\Vert _{k+\left\vert \alpha
\right\vert ,(e_{\ast })}d_{k,e_{\ast }}(\mu _{f},\mu _{g}) \\
&\leq &C\left\Vert a\right\Vert _{d+1}2^{n(k+\left\vert \alpha \right\vert
)}\left\Vert v\right\Vert _{(e_{\ast })}d_{k,e_{\ast }}(\mu _{f},\mu _{g}).
\end{eqnarray*}%
So $\left\Vert \mathcal{H}_{n}^{a}\ast (\partial _{\alpha }(f-g))\right\Vert
_{(e)}\leq C2^{n(k+\left\vert \alpha \right\vert )}d_{k,e_{\ast }}(\mu
_{f},\mu _{g})$ and iv) is proved. $\square $

\smallskip

We are now ready to give the \textquotedblleft balance\textquotedblright .
For $\mu \in \mathcal{M}$ and $\mu _{n}(x)=f_{n}(x)dx,n\in \N$ we recall that
\begin{equation*}
\pi _{q,k,m,e}(\mu ,(\mu _{n})_{n})=\sum_{n=0}^{\infty }2^{n(q+k)}\beta
_{e}(2^{nd})d_{k}(\mu ,\mu _{n})+\sum_{n=0}^{\infty }\frac{1}{2^{2nm}}%
\left\Vert f_{n}\right\Vert _{2m+q,2m,(e)}.
\end{equation*}%
We also set
\begin{equation*}
\widetilde{\pi }_{q,k,m,e}(\mu ,(\mu _{n})_{n})=\sum_{n=0}^{\infty
}2^{n(q+k)}d_{k,e}(\mu ,\mu _{n})+\sum_{n=0}^{\infty }\frac{1}{2^{2nm}}%
\left\Vert f_{n}\right\Vert _{2m+q,2m,(e)}.
\end{equation*}

\begin{proposition}
\label{8}Let $q,k\in \N,m\in \N_{\ast }$ and $e\in \mathcal{E}.$ There exists
a universal constant $C$ (depending on $q,k,m,d$ and $e$) such that for
every $f,f_{n}\in C^{2m+q}(\R^{d}),n\in \N$
\begin{align}
a)\quad & \left\Vert f\right\Vert _{q,(e)}\leq C\pi _{q,k,m,e}(\mu ,(\mu
_{n})_{n}),  \label{Oo10} \\
b)\quad & \left\Vert f\right\Vert _{q,(e)}\leq C\widetilde{\pi }%
_{q,k,m,e}(\mu ,(\mu _{n})_{n}).  \notag
\end{align}%
Here $\mu (x)=f(x)dx$ and $\mu _{n}(x)=f_{n}(x)dx.$


\end{proposition}

\textbf{Proof} \ Let $\alpha $ with $\left\vert \alpha \right\vert \leq q.$
Using Proposition \ref{6}

\begin{equation*}
\partial _{\alpha }f=\sum_{n=1}^{\infty }\mathcal{H}_{n}^{a}\ast \partial
_{\alpha }f=\sum_{n=1}^{\infty }\mathcal{H}_{n}^{a}\ast \partial _{\alpha
}(f-f_{n})+\sum_{n=1}^{\infty }\mathcal{H}_{n}^{a}\ast \partial _{\alpha
}f_{n}
\end{equation*}%
and using (\ref{Oo5}) and (\ref{Oo4})%
\begin{eqnarray*}
\left\Vert \partial _{\alpha }f\right\Vert _{(e)} &\leq &\sum_{n=1}^{\infty
}\left\Vert \mathcal{H}_{n}^{a}\ast \partial _{\alpha }(f-f_{n})\right\Vert
_{(e)}+\sum_{n=1}^{\infty }\left\Vert \mathcal{H}_{n}^{a}\ast \partial
_{\alpha }f_{n}\right\Vert _{(e)} \\
&\leq &C\sum_{n=1}^{\infty }2^{n(\left\vert \alpha \right\vert +k)}\beta
_{e}(2^{nd})d_{k}(\mu _{f},\mu _{f_{n}})+C\sum_{n=1}^{\infty }\frac{1}{%
2^{2nm}}\left\Vert f_{n}\right\Vert _{2m+|\alpha |,2m,(e)}
\end{eqnarray*}%
so (\ref{Oo10}) a) is proved. The same reasoning, using (\ref{Oo6}) gives (%
\ref{Oo10}) b). $\square $

\smallskip

We can now give the

\smallskip

\textbf{Proof of Theorem \ref{2C}.} \textbf{Step 1. Regularization}. For $%
\delta \in (0,1)$ we consider $\gamma _{\delta }$ the density of the centred
Gaussian probability measure with variance $\delta .$ Moreover we consider a
truncation function $\Phi _{\delta }\in C^{\infty }$ such that $1_{B_{\delta
^{-1}}(0)}\leq \Phi _{\delta }\leq 1_{B_{1+\delta ^{-1}}(0)}$ and whose
derivatives of all orders are bounded uniformly w.r.t. $\delta$. Then we
define
\begin{align}
&T_{\delta }\ :\ C^{\infty }\rightarrow C^{\infty }, \quad T_{\delta
}f=(\Phi _{\delta }f)\ast \gamma _{\delta }  \label{Tdelta} \\
&\widetilde{T}_{\delta }\ :\ C^{\infty }\rightarrow C_{c}^{\infty },\quad
\widetilde{T}_{\delta }f=\Phi _{\delta }(f\ast \gamma _{\delta }).  \notag
\end{align}



Moreover, for a measure $\mu \in \mathcal{M}$ we define $T_{\delta }^{\ast
}\mu $ by
\begin{equation*}
\left\langle T_{\delta }^{\ast }\mu ,\phi \right\rangle =\left\langle \mu
,T_{\delta }\phi \right\rangle .
\end{equation*}%
Then $T_{\delta }^{\ast }\mu $ is an absolute continuous measure with
density $p_{T_{\delta }^{\ast }\mu }\in C_{c}^{\infty }$ given by%
\begin{equation*}
p_{T_{\delta }^{\ast }\mu }(y)=\Phi _{\delta }(y)\int \gamma _{\delta
}(x-y)d\mu (x).
\end{equation*}

\textbf{Step 2}. We prove that for every $\mu \in \mathcal{M}$ and $\mu
_{n}(dx)=f_{n}(x)dx,n\in \N$ we have%
\begin{equation}
\pi _{q,k,m,e}(T_{\delta }^{\ast }\mu ,(T_{\delta }^{\ast }\mu
_{n})_{n})\leq C\pi _{q,k,m,e}(\mu ,(\mu _{n})_{n}).  \label{ineg1}
\end{equation}%
Since $\left\Vert T_{\delta }\phi \right\Vert _{k,\infty }\leq C\left\Vert
\phi \right\Vert _{k,\infty }$ one has $d_{k}(T_{\delta }^{\ast }\mu
,T_{\delta }^{\ast }\mu _{n})\leq Cd_{k}(\mu ,\mu _{n}).$

For $\mu _{n}(dx)=f_{n}(x)dx$ we have $p_{T_{\delta }^{\ast }\mu _{n}}(y)=%
\widetilde{T}_{\delta }f_{n}$ Let us now check that
\begin{equation}
\left\Vert \widetilde{T}_{\delta }f_{n}\right\Vert _{2m+q,2m,(e)}\leq
C\left\Vert f_{n}\right\Vert _{2m+q,2m,(e)}.  \label{ineg}
\end{equation}%
For a measurable function $g:\R^{d}\rightarrow \R$ and for $\lambda \geq 0$ we
denote $g_{\lambda }(x)=(1+\left\vert x\right\vert )^{\lambda }g(x).$ Since $%
(1+\left\vert x\right\vert )^{\lambda }\leq (1+\left\vert y\right\vert
)^{\lambda }(1+\left\vert x-y\right\vert )^{\lambda }$ it follows that
\begin{equation*}
(\widetilde{T}_{\delta }g)_{\lambda }(x)=(1+\left\vert x\right\vert
)^{\lambda }\Phi _{\delta }(x)\int_{\R^{d}}\gamma _{\delta }(y)g(x-y)dy\leq
\int_{\R^{d}}\gamma _{\delta ,\lambda }(y)g_{\lambda }(x-y)dy=\gamma _{\delta
,\lambda }\ast g_{\lambda }(x).
\end{equation*}%
Then by (\ref{Oo2}) $\left\Vert (T_{\delta }g)_{\lambda }\right\Vert
_{e}\leq \left\Vert \gamma _{\delta ,\lambda }\ast g_{\lambda }\right\Vert
_{e}\leq C\left\Vert \gamma _{\delta ,\lambda }\right\Vert _{1}\left\Vert
g_{\lambda }\right\Vert _{e}\leq C\left\Vert g_{\lambda }\right\Vert _{e}$.
Using this inequality (with $\lambda =2m)$ for $g=\partial _{\alpha }f_{n}$
we obtain (\ref{ineg}). And (\ref{ineg1}) follows.

\smallskip

\textbf{Step 3}. Let $\mu \in \mathcal{S}_{q,k,m,e}$ so that $\rho
_{q,k,m,e}(\mu )<\infty .$ Using (\ref{Oo10}), a) we have $\left\Vert
T_{\delta }^{\ast }\mu \right\Vert _{W^{q,(e)}}\leq \rho
_{q,k,m,e}(T_{\delta }^{\ast }\mu )$ and moreover, using and (\ref{ineg1})
\begin{equation}
\sup_{\delta \in (0,1)}\left\Vert T_{\delta }^{\ast }\mu \right\Vert
_{W^{q,(e)}}\leq C\sup_{\delta \in (0,1)}\rho _{q,k,m,e}(T_{\delta }^{\ast
}\mu )\leq \rho _{q,k,m,e}(\mu )<\infty .  \label{ineg2}
\end{equation}%
So the family $T_{\delta }^{\ast }\mu ,\delta \in (0,1)$ is bounded in $%
W^{q,e}$ which is a reflexive space. So it is weakly relative compact.
Consequently we may find a sequence $\delta _{n}\rightarrow 0$ such that $%
T_{\delta _{n}}^{\ast }\mu \rightarrow f\in W^{q,e}$ weakly. It is easy to
check that $T_{\delta _{n}}^{\ast }\mu \rightarrow \mu $ weakly so $\mu
(dx)=f(x)dx$ and $f\in W^{q,e}.$ As a consequence of (\ref{ineg2}) we have $%
\left\Vert \mu \right\Vert _{W^{q,(e)}}\leq C\rho _{q,k,m,e}(\mu ).$ $%
\square $

\subsection{Integration by parts formulas}

\label{sect-IBP}

In order to use the criterion presented in the previous section one needs to
have estimates for the densities of the approximating measures $\mu _{n}.$
Sometimes these densities are explicit and then there is no problem. But
sometimes they are not and then the integration by parts machinery is very
useful - we present it in this section. An example related to diffusion
processes is given in Section \ref{sect-diffproc}: if an ellipticity
condition holds then one may exhibit a Gaussian random variable but under H%
\"{o}rmander condition this is no more possible.

We recall that $\mathcal{M}_{+}$ is the set of positive and finite measures
on $\R^{d}$ with the Borel $\sigma $-field and $L_{\mu
}^{p}:=L^{p}(\R^{d},d\mu )$. For $\mu \in \mathcal{M}_{+},m\in \N,p\geq 1$, we
define the Sobolev space $W_{\mu }^{m,p}$ to be the space of the measurable
functions $g:\R^{d}\rightarrow \R$ such that for every multi index $\alpha $
with $\left\vert \alpha \right\vert \leq m$ there exists a function $\theta
_{\alpha }(g)\in L_{\mu }^{p}$ such that
\begin{equation}
\int \partial _{\alpha }f\times gd\mu =(-1)^{\left\vert \alpha \right\vert
}\int f\times \theta _{\alpha }(g)d\mu \quad \forall f\in C_{c}^{\infty
}(\R^{d}).  \label{Not1}
\end{equation}%
We denote
\begin{equation*}
\partial _{\alpha }^{\mu }g=\theta _{\alpha }(g).
\end{equation*}%
Notice that $\partial _{\alpha }^{\mu }$ is not a differential operator.
Indeed it is easy to check that this operator verifies the following
computation rules (see Lemma 9 in \cite{bib:[BCa]}). Let $\phi \in W_{\mu
}^{1,p}$ and $\psi \in C_{b}^{1}(\R^{d}).$ Then $\phi \psi \in W_{\mu }^{1,p}$
and
\begin{equation}
\partial _{i}^{\mu }(\phi \psi )=\psi \partial _{i}^{\mu }\phi +\phi
\partial _{i}\psi .  \label{Not1a}
\end{equation}%
Taking $\phi =1$ we obtain%
\begin{equation}
\partial _{i}^{\mu }\psi =\psi \partial _{i}^{\mu }1+\partial _{i}\psi .
\label{Not1c}
\end{equation}%
We define the Sobolev norm%
\begin{equation*}
\left\Vert g\right\Vert _{W_{\mu }^{m,p}}=\left\Vert g\right\Vert _{L_{\mu
}^{p}}+\sum_{1\leq \left\vert \alpha \right\vert \leq m}\left\Vert \partial
_{\alpha }^{\mu }g\right\Vert _{L_{\mu }^{p}}
\end{equation*}%
and for $p>d$ we denote
\begin{equation}
c_{m,p}(\mu )=\left\Vert 1\right\Vert _{W_{\mu }^{1,p}}^{k_{d,p}}\left\Vert
1\right\Vert _{W_{\mu }^{m,p}}\quad with\quad k_{d,p}=\frac{d-1}{1-d/p}.
\label{Not1d}
\end{equation}

The corresponding definition in terms of random variables is given by means
of the usual integration by parts formulas. On a probability space $(\Omega ,%
\mathcal{F},P)$ we consider a $d$ dimensional $p$ integrable random variable
$F$ and a one dimensional $p$ integrable random variable $G$. Suppose that
for each $\alpha $ with $\left\vert \alpha \right\vert \leq m$ there exists
a $p$ integrable random variable $H_{\alpha }(F,G)$ such that
\begin{equation}
\E(\partial _{\alpha }f(F)G))=(-1)^{\left\vert \alpha \right\vert
}\E(f(F)H_{\alpha }(F,G))\quad \forall f\in C_{c}^{\infty }(\R^{d}).
\label{Not3}
\end{equation}%
If $\mu _{F}$ is the law of $F$ and $g(x)=\E(G\mid F=x)$ then the above
integration by parts formula is equivalent to $g\in W_{\mu _{F}}^{m,p}$ and $%
\theta _{\alpha }(g)=\E(H_{\alpha }(F,G)\mid F=x).$ Indeed, taking
conditional expectations
\begin{eqnarray*}
\int (\partial _{\alpha }f\times g)d\mu _{F} &=&\E(\partial _{\alpha
}f(F)\E(G\mid F))=\E(\partial _{\alpha }f(F)G)) \\
&=&(-1)^{\left\vert \alpha \right\vert }\E(f(F)H_{\alpha
}(F,G))=(-1)^{\left\vert \alpha \right\vert }\E(f(F)\E(H_{\alpha }(F,G)\mid F))
\\
&=&\int (f\times \theta _{\alpha }(g))d\mu _{F}.
\end{eqnarray*}

Moreover, since $\left\vert \E(H_{\alpha }(F,G)\mid F=x)\right\vert ^{p}\leq
\E(\left\vert H_{\alpha }(F,G)\right\vert ^{p}\mid F=x)$ we have%
\begin{equation}
\left\Vert g\right\Vert _{W_{\mu }^{m,p}}^{p}\leq \E\left\vert G\right\vert
^{p}+\sum_{1\leq \left\vert \alpha \right\vert \leq m}\E\left\vert H_{\alpha
}(F,G)\right\vert ^{p}.  \label{Not4}
\end{equation}%
We will express our results in terms of positive measures but, in view of
the above inequality, everything translates immediately in terms of random
variables.

\begin{remark}
An analogous formalism has been introduced by Sigekawa \cite{bib:[Si]} and
by Malliavin and Thalmaier \cite{bib:[MT]} ($\partial _{\alpha }^{\mu }g$
corresponds to the so called \textquotedblleft covering vector
fields\textquotedblright\ from \cite{bib:[MT]}). Here we follow \cite%
{bib:[BCa]}.
\end{remark}

We introduce now the Poisson kernel on $\R^{d}.$ It is given by
\begin{equation*}
Q_{2}(x)=a_{2}\ln \left\vert x\right\vert ,\quad Q_{d}(x)=\frac{a_{d}}{%
\left\vert x\right\vert ^{d-2}},d>2
\end{equation*}%
where $a_{i}$ are some normalization constants. In order to include the one
dimensional case we denote $Q_{1}(x)=x_{+}.$ The kernel $Q_{d}$ is the
fundamental solution to the equation $\triangle Q_{d}=\delta _{0}.$ In \cite%
{bib:[BCa]} Theorem 5 we prove the following estimate. Suppose that $1\in
W_{\mu }^{1,p}$ for some $p>d.$ Then
\begin{equation}
\Theta _{p}(\mu ):=\sup_{x\in \R^{d}}\sum_{i=1}^{d}(\E\int_{\R^{d}}\left\vert
\partial _{i}Q_{d}(x-y)\right\vert ^{\frac{p}{p-1}}d\mu (x))^{(p-1)/p}\leq
C\left\Vert 1\right\Vert _{W_{\mu }^{1,p}}^{k_{d,p}}.  \label{Not7}
\end{equation}%
with $k_{d,p}$ given in (\ref{Not1d}) and $C$ a universal constant.

Following the idea of Malliavin and Thalmaier \cite{bib:[MT]}, in \cite%
{bib:[BCa]} we used the kernel $Q_{d}$ in order to give a representation
theorem for the density of $\mu .$ Consider a function $\psi \in
C_{b}^{\infty }(\R^{d})$ such that $1_{B_{1}(0)}\leq \psi \leq 1_{B_{2}(0)}$
with $B_{r}(0)=\{x:\left\vert x\right\vert <r\}.$ For a fixed $x\in \R^{d}$
we denote $\psi _{x}(y)=\psi (x-y).$ Proposition 9 from \cite{bib:[BCa]}
says that, if $1\in W_{\mu }^{1,p}$ for some $p>d,$ then $\mu (dx)=p_{\mu
}(x)dx$ with
\begin{equation}
p_{\mu }(x)=\sum_{i=1}^{d}\int \partial _{i}Q_{d}(y-x)\partial _{i}^{\mu
}\psi _{x}(y)1_{\{\left\vert x-y\right\vert <2\}}\mu (dy)  \label{Not2}
\end{equation}%
Notice that by (\ref{Not1c}), if $1\in W_{\mu }^{1,p}$ then $\psi _{x}\in
W_{\mu }^{1,p}$ and $\partial _{i}^{\mu }\psi _{x}$ has the support included
in $B_{2}(x).$ So $1_{\{\left\vert x-y\right\vert <2\}}$ appears just to
emphasize this fact.

Suppose now that $1\in W_{\mu }^{m,p}$ for some $m\in \N$ and $p>d.$ In
Theorem 8 from \cite{bib:[BCa]} we prove that $p_{\mu }\in C^{m-1}(\R^{d}).$
And it is easy to see that for a multi index $\alpha $ with $\left\vert
\alpha \right\vert \leq m-1$ one has
\begin{equation*}
\partial _{\alpha }p_{\mu }(x)=(-1)^{\left\vert \alpha \right\vert
}\sum_{i=1}^{d}\int \partial _{i}Q_{d}(y-x)\partial _{\alpha }^{\mu
}\partial _{i}^{\mu }\psi _{x}(y)1_{\{\left\vert x-y\right\vert <2\}}\mu
(dy).
\end{equation*}
We need the following variant of Proposition 9 in \cite{bib:[BCa]}.

\begin{proposition}
\label{1}Suppose that $1\in W_{\mu }^{m,p}$ for some $m\in \N_\ast$ and $p>d.$
Then $\mu (dx)=p_{\mu }(x)dx$ with $p_{\mu }\in C^{m-1}(\R^{d})$ and for
every multi index $\alpha $ with $\left\vert \alpha \right\vert \leq m-1$
and $k\in \N$%
\begin{equation}
\left\vert \partial _{\alpha }p_{\mu }(x)\right\vert \leq C\times
c_{m,2p}(\mu )m_{k}^{\frac{1}{2}}(\mu )u_{k/2}(x)  \label{Not9a}
\end{equation}%
with $u_{k}(x)=(1+\left\vert x\right\vert )^{-k}$ and
\begin{equation*}
m_{k}(\mu )=\int (1+\left\vert x\right\vert )^{k}d\mu (x).
\end{equation*}
\end{proposition}

\textbf{Proof}. We write%
\begin{equation*}
\left\vert \partial _{\alpha }p_{\mu }(x)\right\vert \leq
\sum_{i=1}^{d}I_{i}\times J_{i}\leq
(\sum_{i=1}^{d}I_{i}^{2})^{1/2}(\sum_{i=1}^{d}J_{i}^{2})^{1/2}
\end{equation*}%
with
\begin{equation*}
I_{i}=(\int \left\vert \partial _{i}Q_{d}(y-x)\right\vert ^{\frac{p}{p-1}%
}\mu (dy))^{\frac{p-1}{p}}\quad and\quad J_{i}=(\int \left\vert \partial
_{\alpha }^{\mu }\partial _{i}^{\mu }\psi _{x}(y)\right\vert
^{p}1_{\{\left\vert x-y\right\vert <2\}}\mu (dy))^{\frac{1}{p}}.
\end{equation*}%
We use (\ref{Not7}) and we obtain
\begin{equation*}
(\sum_{i=1}^{d}I_{i}^{2})^{1/2}\leq C\left\Vert 1\right\Vert _{W_{\mu
}^{1,p}}^{k_{d,p}}.
\end{equation*}%
Moreover using Schwarz inequality
\begin{equation*}
J_{i}\leq (\int \left\vert \partial _{\alpha }^{\mu }\partial _{i}^{\mu
}\psi _{x}(y)\right\vert ^{2p}\mu (dy))^{\frac{1}{2p}}\mu (B_{2}(x))^{\frac{1%
}{2}}.
\end{equation*}%
Using repeatedly (\ref{Not1c}) one can check that $\left\Vert \psi
_{x}\right\Vert _{W_{\mu }^{m,2p}}\leq C_{\psi }\left\Vert 1\right\Vert
_{W_{\mu }^{m,2p}}$ where $C_{\psi }$ is a universal constant which depends
on the derivatives of $\psi $ up to order $m.$ Take now a random variable $%
\xi $ of law $\mu .$ We may assume without loss of generality that $%
\left\vert x\right\vert \geq 4$ and then, using Chebyshev inequality%
\begin{equation*}
\mu (B_{2}(x))=\P(\xi \in B_{2}(x))\leq \P(\left\vert \xi \right\vert \geq
\frac{\left\vert x\right\vert }{2})\leq \frac{C}{\left\vert x\right\vert ^{k}%
}m_{k}(\mu ).
\end{equation*}

So $J_{i}\leq C\left\Vert \mu \right\Vert _{m,2p}m_{k}^{1/2}(\mu
)u_{k/2}(x). $ We conclude that
\begin{equation*}
\left\vert \partial _{\alpha }p_{\mu }(x)\right\vert \leq C\left\Vert
1\right\Vert _{W_{\mu }^{1,p}}^{k_{d,p}}\left\Vert 1\right\Vert _{W_{\mu
}^{m,2p}}m_{k}^{1/2}(\mu )u_{k/2}(x)=c_{m,2p}(\mu )m_{k}^{1/2}(\mu
)u_{k/2}(x)
\end{equation*}%
and consequently (\ref{Not9a}) is proved. $\square $

\smallskip

By using Remark \ref{U}, as an immediate consequence of (\ref{Not9a}) (with $%
p=d+1)$ we obtain%
\begin{equation}
\left\Vert p_{\mu }\right\Vert _{2m+q,2m,(e)}\leq C_{d,e}c_{2m+q,2(d+1)}(\mu
)m_{2(d+1+m)}^{1/2}(\mu ).  \label{Not20}
\end{equation}

We can now re-formulate the results from Section \ref{sect-results}. We
define
\begin{equation*}
\widetilde{M}_{m,q,e}(R)=\{\mu \in \mathcal{M}_{a}(\R^{d}):c_{2m+q,2(d+1)}(%
\mu )m_{2(d+1+m)}^{1/2}(\mu )\leq R\}
\end{equation*}%
and by (\ref{Not9a}) we obtain $\widetilde{M}_{m,q,e}(R)\subset
M_{m,q,e}(C\times R)$ where $C$ is a universal constant. Then we consider
the following

\smallskip

\textbf{Hypothesis $\widetilde{H}_{q}(k,m,e)$}. \emph{For $q,k\in \N$ and $%
m\in \N_{\ast }$ there exists $a>1$ such that%
\begin{equation*}
\overline{\lim }_{R\rightarrow \infty }\frac{L_{a}(R)^{1+\frac{k+q}{2m}%
}\beta _{e}(L_{a}(R)^{\frac{d}{2m}})}{R}d_{k}(\mu ,\widetilde{M}%
_{m,q,e}(R))<\infty .
\end{equation*}%
}

\medskip

We define
\begin{equation*}
\widetilde{B}_{q}(k,m,e)=\{\mu \in \mathcal{M}(\R^{d}):\text{ }\widetilde{H}%
_{q}(k,m,e)\text{ holds for }\mu \}.
\end{equation*}

Our criterion can be stated as follows.

\begin{theorem}
\label{2F}Let $q,k\in \N$, $m\in \N_{\ast }$ and $e\in \mathcal{E}$. Then $%
\widetilde{B}_{q}(k,m,e)\subset B_{q}(k,m,e)\subset \mathcal{S}%
_{q,k,m,e}\subset W^{q,e}.$
\end{theorem}

\section{Interpolation spaces}

\label{sect-interp}

We give a variant of the results from the previous section in terms of
interpolation spaces. More precisely we will prove that the space $\mathcal{S%
}_{q,k,m,e}$ considered in Section \ref{sect-results} (\ref{O13}) is an
interpolation space.

To begin we recall the framework of interpolation spaces. We are given two
Banach spaces $(X,\left\Vert \circ \right\Vert _{X}),(Y,\left\Vert \circ
\right\Vert _{Y})$ with $X\subset Y$ (with continuous embedding). We denote $%
\mathcal{L}(X,X)$ the space of the linear bounded operators from $X$ into
itself and we denote by $\left\Vert L\right\Vert _{X,X}$ the operator norm.
A Banach space $(W,\left\Vert \circ \right\Vert _{W})$ such that $X\subset
W\subset Y$ is called an interpolation space for $X$ and $Y$ if $\mathcal{L}%
(X,X)\cap \mathcal{L}(Y,Y)\subset \mathcal{L}(W,W).$ Let $\gamma \in (0,1).$
If there exists a constant $C$ such that $\left\Vert L\right\Vert _{W,W}\leq
C\left\Vert L\right\Vert _{X,X}^{\gamma }\left\Vert L\right\Vert
_{Y,Y}^{1-\gamma }$ for every $L\in \mathcal{L}(X,X)\cap \mathcal{L}(Y,Y)$
then $W$ is an interpolation space of order $\gamma .$ And if one may take $%
C=1$ then $W$ is an exact interpolation space of order $\gamma .$ There are
several methods for constructing interpolation spaces. We focus here on the
so called $K$-method. For $y\in Y$ and $t>0$ one defines $K(y,t)=\inf_{x\in
X}(\left\Vert y-x\right\Vert _{Y}+t\left\Vert x\right\Vert _{X})$ and
\begin{equation*}
\left\Vert y\right\Vert _{\gamma }=\int_{0}^{\infty }t^{-\gamma }K(y,t)\frac{%
dt}{t},\qquad (X,Y)_{\gamma }=\{y\in Y:\left\Vert y\right\Vert _{\gamma
}<\infty \}.
\end{equation*}%
Then one proves that $(X,Y)_{\gamma }$ is an exact interpolation space of
order $\gamma .$

In the theory of interpolation spaces one considers a more general
situation: one does not require that $X\subset Y$ but instead one works with
$X\cap Y$ and $X+Y.$ Moreover one considers a whole family of norms given by
$\left\Vert y\right\Vert _{\gamma ,q}^{q}=\int_{0}^{\infty }(t^{-\gamma
}K(y,t))^{q}\frac{dt}{t}$ with $\infty >q\geq 1.$ In our framework we are
concerned with $X\subset Y$ and $q=1$ so we leave out the general setting.
But we will need a variant of the norm $\left\Vert y\right\Vert _{\gamma }:$
for $\gamma \in (0,1)$ and $b\geq 0$ we define
\begin{equation*}
\left\Vert y\right\Vert _{\gamma ,b}=\int_{0}^{\infty }t^{-\gamma
}\left\vert \ln t\right\vert ^{b}K(y,t)\frac{dt}{t}\qquad and\qquad
\left\vert y\right\vert _{\gamma ,b}=\int_{0}^{1}t^{-\gamma }\left\vert \ln
t\right\vert ^{b}K(y,t)\frac{dt}{t}.
\end{equation*}%
Notice that $K(y,t)\leq \left\Vert y\right\Vert _{Y}$ so if $\gamma >0$ then%
\begin{equation*}
\left\vert y\right\vert _{\gamma ,b}\leq \left\Vert y\right\Vert _{\gamma
,b}\leq \left\vert y\right\vert _{\gamma ,b}+C_{\gamma ,b}\Vert y\Vert _{Y}.
\end{equation*}%
Therefore we may work with one or another of these two norms. But if $\gamma
=0$ then $\int_{1}^{\infty }\left\vert \ln t\right\vert ^{b}K(t,y)\frac{dt}{t%
}$ may be infinite. This is why we prefer to work with $\left\vert
y\right\vert _{\gamma ,b}.$ We define%
\begin{equation*}
K_{\gamma ,b}(X,Y)=\{y\in Y:\left\vert y\right\vert _{\gamma ,b}<\infty \}.
\end{equation*}%
In the case $\gamma >0$ we have $K_{\gamma ,b}(X,Y)=\{y\in Y:\left\Vert
y\right\Vert _{\gamma ,b}<\infty \}$ and moreover if $b=0$ we have $%
\left\Vert y\right\Vert _{\gamma }=\left\Vert y\right\Vert _{\gamma ,0}$ and
$(X,Y)_{\gamma }=K_{\gamma ,0}(X,Y)$.

\smallskip

We introduce now another space which is the analogous of the space $\mathcal{%
S}_{q,k,m,e}$ defined in (\ref{O13}). Let $a,\theta \geq 0$ and $m\in \N_{\ast }.$ For $y\in Y$ and for a sequence $x_{n}\in X,n\in \N$ we define
\begin{equation}
\pi _{\theta ,m,a}(y,(x_{n})_{n})=\sum_{n=1}^{\infty }2^{n\theta
}n^{a}\left\Vert y-x_{n}\right\Vert _{Y}+\frac{1}{2^{2nm}}\left\Vert
x_{n}\right\Vert _{X}.  \label{Int1}
\end{equation}%
Moreover we define
\begin{equation*}
\rho _{\theta ,m,a}^{X,Y}(y)=\inf \pi _{\theta ,m,a}(y,(x_{n})_{n})
\end{equation*}%
with the infimum taken over all the sequences $x_{n}\in X,n\in \N.$ We define
\begin{equation}
S_{\theta ,m,a}(X,Y)=\{y\in Y:\rho _{\theta ,m,a}^{X,Y}(y)<\infty \}.
\label{Int2}
\end{equation}%
It is easy to check that $\rho^{X,Y} _{\theta ,m,a}$ is a norm on $S_{\theta
,m,a}(X,Y).$

Finally we define the space which corresponds to the balance: for $\alpha
,\beta \geq 0$
\begin{equation*}
B_{\alpha ,\beta }(X,Y)=\{y\in Y:\overline{\lim }_{R\rightarrow \infty
}R^{\alpha }(\ln R)^{\beta }d_{Y}(y,B_{X}(R))<\infty \}.
\end{equation*}%
In the appendix we prove the following relation between these spaces (see
Proposition \ref{NORM} and Proposition \ref{balance}):

\begin{proposition}
\label{Spaces}Given $\theta ,a\geq 0$ and $m\in \N_{\ast }$, one has
\begin{equation}
B_{\alpha ,\beta }(X,Y)\subset S_{\theta ,m,a}(X,Y)=K_{\gamma ,b}(X,Y)
\label{Int41}
\end{equation}%
with%
\begin{equation*}
\gamma =\frac{\theta }{2m+\theta },b=\frac{2ma}{2m+\theta }\qquad and\qquad
\alpha =\frac{\theta }{2m},\beta =2+a+\frac{\theta }{m}.
\end{equation*}%
In particular, taking $a=0$ we obtain%
\begin{equation*}
B_{\alpha ,\beta }(X,Y)\subset (X,Y)_{\gamma }\qquad with\qquad \alpha =%
\frac{\gamma }{1-\gamma },\beta =\frac{2}{1-\gamma }.
\end{equation*}
\end{proposition}

\subsection{Distribution spaces}

We precise now the particular spaces we will work with. We denote by $%
W^{k,\infty }$ the closure of the $C_{b}^{\infty }$ functions with respect
to $\left\Vert \circ \right\Vert _{k,\infty }$ and we consider the dual
space $W_{\ast }^{k,\infty }.$ So $W_{\ast }^{k,\infty }$ is the space of
the bounded linear functionals $u:W^{k,\infty }\rightarrow \R$ with the norm
\begin{equation*}
\left\Vert u\right\Vert _{W_{\ast }^{k,\infty }}=\sup \{\left\vert
\left\langle u,f\right\rangle \right\vert :f\in W^{k,\infty },\left\Vert
f\right\Vert _{k,\infty }\leq 1\}.
\end{equation*}

Let $k,q,m\in \N,m>\frac{d}{2}$ and $e\in \mathcal{E}_{\alpha ,\gamma }$ for
some $\alpha ,\gamma \geq 0$ (see (\ref{RR}))$.$\ We will work with the
spaces%
\begin{equation}
X=W^{2m+q,2m,e}\qquad and\qquad Y=W_{\ast }^{k,\infty }.  \label{S}
\end{equation}%
It is easy to check that, if $l>d,$ then $W^{0,l,e}\subset L^{1}.$ Indeed if
$u_{l}(x)=(1+\left\vert x\right\vert )^{-l}$ and $v(x)=(1+\left\vert
x\right\vert )^{l}f(x)$ then using H\"{o}lder inequality $\int \left\vert
f(x)\right\vert dx=\int \left\vert u_{l}(x)v(x)\right\vert dx\leq
2\left\Vert v\right\Vert _{(e)}\left\Vert u_{l}\right\Vert _{(e_{\ast
})}\leq 2\left\Vert f\right\Vert _{0,l,(e)}\left\Vert u_{l}\right\Vert
_{(e_{\ast })}$ and by Remark \ref{U} this is finite. So we may embed $%
X=W^{2m+q,2m,e}$ in $Y=W_{\ast }^{k,\infty }$ in the sense that $f\in
X\rightarrow \widetilde{f}\in Y$ with $\langle \widetilde{f},g\rangle =\int
f(x)g(x)dx.$ Moreover for a measure $\mu \in \mathcal{M}$ we denote by $%
\widetilde{\mu }$ the linear functional $\left\langle \widetilde{\mu }%
,g\right\rangle =\int g(x)d\mu (x).$ Since $g\in W^{k,\infty }$ is bounded $%
\left\langle \widetilde{\mu },g\right\rangle $ is well defined. So one has $%
\widetilde{\mu }\in W_{\ast }^{k\infty ,}$. The distance introduced in (\ref%
{O6}) becomes
\begin{eqnarray*}
d_{k}(\mu ,\nu ) &=&\sup \{\left\vert \int \phi d\mu -\int \phi d\nu
\right\vert :\phi \in C^{k}(\R^{d}),\left\Vert \phi \right\Vert _{k,\infty
}\leq 1\} \\
&=&\sup \{\left\vert \left\langle \widetilde{\mu }-\widetilde{\nu },\phi
\right\rangle \right\vert :\phi \in C^{k}(\R^{d}),\left\Vert \phi \right\Vert
_{k,\infty }\leq 1\}=\left\Vert \widetilde{\mu }-\widetilde{\nu }\right\Vert
_{W_{\ast }^{k,\infty }}.
\end{eqnarray*}%
Recall $\rho _{q,k,m,e}$ defined in (\ref{011}) and take $\theta =q+k+\alpha
d,a=\gamma .$ Given $e\in \mathcal{E}_{\alpha ,\gamma }$ one may find a
universal $C$ such that for every $\mu \in \mathcal{M}$ one has
\begin{equation}
\frac{1}{C}\rho _{q,k,m,e}(\mu )\leq \rho _{\theta ,m,a}^{X,Y}(\widetilde{%
\mu })\leq C\rho _{q,k,m,e}(\mu ).  \label{ineg4}
\end{equation}%
We have the following analogous of Theorem \ref{2C}:

\begin{proposition}
Let $k,q,m\in \N,m>\frac{d}{2}$ and $e\in \mathcal{E}_{\alpha ,\gamma }$. We
take $\theta =q+k+\alpha d$\ and $X=W^{2m+q,2m,e},Y=W_{\ast }^{k,\infty }.$
Then,
\begin{equation*}
S_{\theta ,m,\gamma }(X,Y)\subset W^{q,e}.
\end{equation*}
\end{proposition}

The proof is identical with the proof of Theorem \ref{2C} (we have not used
there the fact that $\mu $ is a measure, but only the fact that it is a
distribution) so we skip it.

The main result in this section is the following:

\begin{theorem}
\label{th3.3} Let $k,q,m\in \N$ and $e\in \mathcal{E}_{\alpha ,\gamma }$ for
some $\alpha ,\gamma \geq 0.$ We take $\theta =q+k+\alpha d$ and we denote $%
X=W^{2m+q,2m,e}\subset W_{\ast }^{k,\infty }=Y.$ Then,
\begin{align*}
& W^{q+1,2m,e}\subset \mathcal{S}_{q,k,m,e}=S_{\theta ,m,\gamma
}(X,Y)=K_{\rho ,b}(X,Y)\subset W^{q,e} \\
& with\quad \rho =\frac{\theta }{2m+\theta },b=\frac{2m\gamma }{2m+\theta },
\end{align*}%
the first inclusion holding for $m>d/2$, and
\begin{equation*}
B_{u,v}(X,Y)\subset S_{\theta ,m,\gamma }(X,Y)\qquad u=\frac{\theta }{2m}%
,v=2+\gamma +\frac{\theta }{m}.
\end{equation*}
\end{theorem}

\textbf{Proof}. As a consequence of the previous proposition the elements of
$S_{\theta ,m,\gamma }(X,Y)$\ are functions.\ So using (\ref{ineg4}) we
obtain $\mathcal{S}_{q,k,m,e}=S_{\theta ,m,\gamma }(X,Y).$ The equality $%
S_{\theta ,m,\gamma }(X,Y)=K_{\rho ,b}(X,Y)$ and the inclusion $%
B_{u,v}(X,Y)\subset S_{\theta ,m,\gamma }(X,Y)$ are proved in (\ref{Int41}).
The inclusions $W^{q+1,2m,e}\subset S_{\theta ,m,\gamma }(X,Y)=K_{\rho
,b}(X,Y)\subset W^{q,e}$ follow from Theorem \ref{2E} and the fact that $%
S_{\theta ,m,\gamma }(X,Y)=\mathcal{S}_{q,k,m,e}.$ $\square $

\smallskip

We come now back to the examples given in the first section.

\smallskip

\textbf{Example 1}. \emph{We take $e_{p}(t)=\left\vert t\right\vert ^{p}$ so
that $\beta _{e_{p}}(t)=\left\vert t\right\vert ^{1/p_{\ast }}.$ It follows
that $e_{p}\in \mathcal{E}_{\alpha ,\gamma }$ with $\alpha =\frac{1}{p_{\ast
}}$ and $\gamma =0.$ The Sobolev spaces associated to $e_{p}$ are the
standard Sobolev spaces so we denote $W^{q,p}=W^{q,e_{p}}$ and so on. We
also have $b=0$ so $K_{\rho ,b}(X,Y)=(X,Y)_{\rho }$ which is the standard
interpolation space. Moreover $\theta =q+k+\alpha d=q+k+d/p_{\ast }\rho
=\theta /(2m+\theta )$ and we obtain%
\begin{align*}
& W^{q+1,2m,p}\subset \mathcal{S}_{q,k,m,e_{p}}=(W^{2m+q,2m,p},W_{\ast
}^{k,\infty })_{\rho }\subset W^{q,p}, \\
& \text{with }\rho =\frac{q+k+d/p_{\ast }}{2m+q+k+d/p_{\ast }}.
\end{align*}%
}

\smallskip

\textbf{Example 2}. \emph{We take $e_{\log }(t)=(1+\left\vert t\right\vert
)\ln (1+\left\vert t\right\vert ).$ So $e_{\log }\in \mathcal{E}_{0,1}$ $.$
Now we work with%
\begin{equation*}
K_{\rho ,b}(X,Y)=\{y\in Y:\int_{0}^{1}\frac{\left\vert \ln t\right\vert ^{b}%
}{t^{\rho }}K(y,t)\frac{dt}{t}<\infty \}.
\end{equation*}%
This is no more a standard interpolation space. The above theorem gives%
\begin{align*}
&W^{q+1,2m,e_{\log}} \subset \mathcal{S}_{q,k,m,e_{\log }}=K_{\rho
,b}(W^{2m+q,2m,e_{\log }},W_{\ast }^{k,\infty })\subset W^{q,e_{\log }}, \\
&\text{with } \rho =\frac{q+k}{2m+q+k} \text{ and } b=\frac{2m}{2m+q+k}.
\end{align*}
}

\subsection{Negative Sobolev spaces}

Usually in interpolation theory one discusses about $L^{p}$ spaces and the
case of $L^{\infty }$ treated above appears as a limit case which is more
delicate. In this section we investigate the results that one may obtain for
$L^{p}$ using the machinery based on Hermite expansions (presented in
Section \ref{sect-Hermite}). For simplicity we leave out the Orlicz norms
and we restrict ourself to $L^{p}$ spaces. So the spaces we work with are $%
X=W^{2m+q,2m,p}$ and $Y=W^{-k,p_{\ast }}=(W^{k,p_{\ast }})^{\ast }$ the dual
of the standard Sobolev space $W^{k,p_{\ast }}.$ The result in this section
is the following:

\begin{theorem}
Let $k,q,m\in \N.$ We denote $X=W^{2m+q,2m,p}\subset W^{-k,p_{\ast }}=Y.$ Then%
\begin{equation*}
W^{q+1,2m,e}\subset (W^{2m+q,2m,p},W^{-k,p_{\ast }})_{\rho }\subset
W^{q,p}\qquad with\qquad \rho =\frac{q+k}{2m+q+k},
\end{equation*}%
first inclusion holding only for $p_{\ast }>\frac{d(2m-1)}{2m+k+q}$ and $%
m>d/2$.
\end{theorem}

\textbf{Proof}. The proof is the same as for Theorem \ref{2C} so we sketch
it only (we use the notation introduced in that proof).

Let us first prove that $(X,Y)_{\rho }\subset W^{q,p}.$ By (\ref{Int41}) we
have $(X,Y)_{\rho }=K_{\rho ,0}=S_{\theta ,m,0}(X,Y)$ with $\theta =2m\rho
/(1-\rho ).$ Let $u\in (X,Y)_{\rho }.$ By the very definition of $S_{\theta
,m,0}(X,Y)$ we may find a sequence $v_{n}\in X$ such that
\begin{equation*}
\pi _{\theta ,m,0}(u,(v_{n})_{n})=\sum_{n=1}^{\infty }2^{n\theta }\left\Vert
u-v_{n}\right\Vert _{Y}+\frac{1}{2^{2nm}}\left\Vert v_{n}\right\Vert
_{X}<\infty .
\end{equation*}

For $u\in Y=(W^{k,p_{\ast }})^{\ast }$\ we define $T_{\delta }^{\ast }u$ by $%
\left\langle T_{\delta }^{\ast }u,\phi \right\rangle =\left\langle
u,T_{\delta }\phi \right\rangle $ - the operators $T_{\delta },T_{\delta
}^{\ast },\widetilde{T}_{\delta }$ are defined in the Step 1 of the proof of
Theorem \ref{2C}), see (\ref{Tdelta}). The functional $T_{\delta }^{\ast }u$
is associated to the function $f(y):=\Phi _{\delta }(y)\left\langle u,\gamma
_{\delta }(\cdot -y)\right\rangle $ in the sense that $\left\langle
T_{\delta }^{\ast }u,\phi \right\rangle =\int f(y)\phi (y)dy.$

Since $\left\Vert T_{\delta }\phi \right\Vert _{W^{k,p_{\ast }}}\leq
\left\Vert \phi \right\Vert _{W^{k,p_{\ast }}}$ we obtain $\left\Vert
T_{\delta }^{\ast }u-T_{\delta }^{\ast }v_{n}\right\Vert _{Y}\leq \left\Vert
u-v_{n}\right\Vert _{Y}.$ Notice that the notation $T_{\delta }^{\ast }v_{n}$
is introduced for the functional $\phi \rightarrow \int v_{n}(y)\phi (y)dy.$
So $\left\langle T_{\delta }^{\ast }v_{n},\phi \right\rangle =\left\langle
f_{n},\phi \right\rangle $ with $f_{n}(y)=\Phi _{\delta }(y)\left\langle
v_{n},\gamma _{\delta }(\cdot -y)\right\rangle =\Phi _{\delta }(y)v_{n}\ast
\gamma _{\delta }(y)=\widetilde{T}_{\delta }v_{n}.$ We have already proved
in (\ref{ineg}) (take $f_{n}=v_{n})$ that $\left\Vert \widetilde{T}_{\delta
}v_{n}\right\Vert _{X}\leq \left\Vert v_{n}\right\Vert _{X}.$ We conclude
that%
\begin{equation*}
\pi _{\theta ,m,0}(T_{\delta }^{\ast }u,(T_{\delta }^{\ast }v_{n})_{n})\leq
\pi _{\theta ,m,0}(u,(v_{n})_{n}).
\end{equation*}

Using (\ref{Oo10}), b) we obtain
\begin{equation*}
\left\Vert T_{\delta }^{\ast }u\right\Vert _{W^{q,p}}\leq C\pi _{\theta
,m,0}(u,(v_{n})_{n})
\end{equation*}%
so the family $T_{\delta }^{\ast }u,\delta >0$ is bounded and consequently
relatively compact in $W^{q,p}.$ Let $\widetilde{u}\in W^{q,p}$ be a limit
point. Since $\left\langle T_{\delta }^{\ast }u,\phi \right\rangle
\rightarrow \left\langle u,\phi \right\rangle $ we have $\left\langle u,\phi
\right\rangle =\left\langle \widetilde{u},\phi \right\rangle .$

The proof of the inclusion $W^{q+1,2m, e}\subset
(W^{2m+q,2m,p},W^{-k,p_{\ast }})_{\rho }$ is analogous to that of (\ref{K})
so we skip it. $\square $

\begin{remark}
Let us compare this result with the one in Example 1 given before. We have
replaced the dual of $W^{k,\infty }$ by the dual of $W^{k,p_{\ast }}$ and
then the interpolation index is smaller and does no more depend on the
dimension (because $d/p_{\ast }$ does no more appear).
\end{remark}

\section{Diffusion processes}

\label{sect-diffproc}

We consider the SDE%
\begin{equation}
X_{t}=x+\sum_{j=1}^{N}\int_{0}^{t}\sigma
_{j}(s,X_{s})dW_{s}^{j}+\int_{0}^{t}b(s,X_{s})ds  \label{D1}
\end{equation}%
with $\sigma _{j},b:\R_{+}\times \R^{d}\rightarrow \R^{d},j=1,...,N$ measurable
functions$\ $and $W$ a standard $N$-dimensional Brownian motion. We want to keep weak
hypotheses on the coefficients so we do not know that the above SDE has
a unique solution. So, we just consider a continuous and adapted process $%
X_{t}$, $t\geq 0$, which verifies (\ref{D1}). We assume that the coefficients
have linear growth: for every $T>0$ there exists a constant $C_{T}$ such
that
\begin{equation}
\left\vert b(t,x)\right\vert +\sum_{j=1}^{N}\left\vert \sigma
_{j}(t,x)\right\vert \leq C_{T}(1+\left\vert x\right\vert )\quad \forall
(t,x)\in \lbrack 0,T]\times \R^{d}.  \label{D2}
\end{equation}%
Moreover for an open domain $D\subset \R^{d}$ we define $C_{\log
}([0,T]\times D)$ to be the set of functions $f:[0,T]\times D\rightarrow
\R^{d}$ for which there exist $C,h>0$ such that for every $(t,x),(s,y)\in
\lbrack 0,T]\times D$
\begin{equation}
\left\vert f(t,x)-f(s,y)\right\vert \leq C(\left\vert \ln \left\vert
x-y\right\vert \right\vert ^{-(2+h)}+\left\vert \ln \left\vert
t-s\right\vert \right\vert ^{-(2+h)}).  \label{D3}
\end{equation}%
Finally we say that a probability measure $\mu $ on $\R^d$ has a local density $p_{\mu
}$ on $D$ if $\mu (A)=\int_{A}p_{\mu }(x)dx$ for every Borel set $A\subset D.$ We consider
now an open domain $\Gamma \subset \R^{d}$ and we define
\begin{equation*}
\tau =\inf \{t:X_{t}\notin \Gamma \}.
\end{equation*}

\subsection{A local ellipticity condition}

We first study the case in which the diffusion coefficients $\sigma_j$, $j=1,\ldots,N$, satisfy a local ellipticity condition. More precisely, we prove that
\begin{theorem}
\label{9}
Let the coefficients $b$ and $\sigma _{j}$, $j=1,\ldots,N$, be measurable and with
linear growth. For $y_{0}\in \Gamma$, $T>0$ and $0<r<\frac{1}{2}%
d(y_{0},\Gamma ^{c})$, suppose that $\sigma _{j}\in C_{\log
}([0,T]\times B_{r}(y_{0})),j=1,...,N$ and $\sigma \sigma ^{\ast }(y_{0})>0.$
Let $X_{t}$\ be any continuous and adapted solution to the SDE (\ref{D1}).
Then the law of $X_{T\wedge \tau}$ has a local density on $%
B_{r/4}(y_{0}).$
\end{theorem}

We will need the following approximation result. We fix $\delta >0$ and we
construct $X_{t}^{\delta }$ by%
\begin{equation*}
X_{t}^{\delta }=X_{(T-\delta )\wedge \tau }+\sum_{j=1}^{N}\sigma (T-\delta
,X_{T-\delta })(W_{t}^{j}-W_{T-\delta }^{j}),\quad t\in[T-\delta,T].
\end{equation*}%
We denote by $\mu $ the law of $X_{T\wedge \tau}$ and by $\mu
_{\delta }$ the law of $X_{T}^{\delta }.$ Moreover we consider a truncation
function $\psi \in C^{\infty }$ such that $1_{B_{r/2}(y_{0})}\leq \psi \leq
1_{B_{r}(y_{0})}$ and such that for every multi index $\alpha $ one has $%
\left\Vert \partial ^{\alpha }\psi \right\Vert _{\infty }\leq C_{\alpha
}r^{-\left\vert \alpha \right\vert }$ where $C_{\alpha }$ depends on $\alpha
$ but not on $r.$ We define $\nu =\psi \mu $ and $\nu _{\delta }=\psi \mu
_{\delta }.$

\begin{lemma}\label{10}
Under the hypothesis of the Theorem \ref{9} one has
\begin{equation*}
d_{1}(\nu ,\nu _{\delta })\leq \frac{C_*}{r}\delta ^{1/2}(\ln \delta )^{-(2+h)}
\end{equation*}%
where $C_*>0$ depends on the constants in (\ref{D2}) and (\ref%
{D3}) associated to $\sigma _{j}$ and $b$, $h$ being the constant from (\ref{D3}) related to $\sigma _{j}$ and $b$.
\end{lemma}

\textbf{Proof}. Since $X_{\tau }\in \Gamma ^{c}$ and $r\leq \frac{1}{2}%
d(y_{0},\Gamma ^{c})$ we have $\psi (X_{\tau })=0$ so $\psi (X_{T\wedge \tau
})=\psi (X_{T})1_{\{\tau >T\}}.$ It follows that for any $\phi \in
C_{b}^{1}(\R^{d})$
\begin{equation*}
\int \phi d\nu =\E(\phi \psi (X_{T\wedge \tau })=\E(\phi \psi (X_{T})1_{\{\tau
>T\}})=R_{\delta }+I
\end{equation*}%
with
\begin{equation*}
R_{\delta }=\E((\phi \psi (X_{T})-\phi \psi (X_{T}^{\delta }))1_{\{\tau
>T\}})\quad and\quad I=\E(\phi \psi (X_{T}^{\delta })1_{\{\tau >T\}}).
\end{equation*}%
We write%
\begin{equation*}
I=\int \phi d\nu _{\delta }-J\quad with\quad J=\E(\phi \psi (X_{T}^{\delta
})1_{\{\tau \leq T\}}).
\end{equation*}%
So%
\begin{equation*}
\left\vert \int \phi d\nu -\int \phi d\nu _{\delta }\right\vert \leq
\left\vert J\right\vert +\left\vert R_{\delta }\right\vert .
\end{equation*}

We estimate $\left\vert J\right\vert .$\ Let $\eta =\inf \{t>T-\delta
:\left\vert X_{t}-X_{(T-\delta )\wedge \tau }\right\vert >\frac{r}{4}%
\} $ and $\eta _{\delta }=\inf \{t>T-\delta :\left\vert X_{t}^{\delta
}-X_{(T-\delta )\wedge \tau }\right\vert >\frac{r}{4}\}.$ Using
standard arguments one checks that
\begin{equation*}
\P(\eta <T)+\P(\eta _{\delta }<T)\leq C\delta .
\end{equation*}%
Suppose that we are on the set $\{\eta >T\}\cap \{\eta _{\delta }>T\}.$ If $%
\psi (X_{T}^{\delta })\neq 0$ then $\left\vert X_{T}^{\delta
}-y_{0}\right\vert \leq \frac{r}{2}$ and since $\eta _{\delta }>T$ it
follows that $\left\vert X_{(T-\delta )\wedge \tau }-y_{0}\right\vert \leq
\frac{3r}{4}<r.$ Since $\left\vert X_{\tau }-y_{0}\right\vert >r$ this
implies that $\tau >T-\delta .$ And we have $\left\vert X_{T-\delta
}-y_{0}\right\vert \leq \frac{3r}{4}.$ And since $\eta >T$ it follows also
that $\left\vert X_{t}-y_{0}\right\vert \leq r$ for every $T-\delta \leq
t\leq T.$ Then $X_{t}\in \Gamma $ and so $\tau >T.$ We conclude that $\{\psi
(X_{T}^{\delta })\neq 0\}\cap \{\tau \leq T\}\subset \{\eta <T\}\cup \{\eta
_{\delta }<T\}.$ So that
\begin{equation*}
\left\vert J\right\vert \leq \left\Vert \phi \psi \right\Vert _{\infty
}(\P(\eta <T)+\P(\eta _{\delta }<T))\leq C\left\Vert \phi \right\Vert _{\infty
}\delta .
\end{equation*}%
We estimate now%
\begin{eqnarray*}
\left\vert R_{\delta }\right\vert &\leq &\left\Vert \nabla (\phi \psi
)\right\Vert _{\infty }\E(\left\vert X_{T}-X_{T}^{\delta }\right\vert
1_{\{\tau >T\}}1_{\{X_{T}\in B_{r}(y_{0})\}}) \\
&\leq &\frac{C}{r}(\left\Vert \phi \right\Vert _{\infty }+\left\Vert \nabla
\phi \right\Vert _{\infty })\E(\left\vert X_{T}-X_{T}^{\delta }\right\vert
1_{\{\eta >T\}\cap \{\tau >T\}\cap \{X_{T}\in B_{r}(y_{0})\}})+C\left\Vert
\phi \right\Vert _{\infty }\delta .
\end{eqnarray*}%
On the set $\{\eta >T\}\cap \{\tau >T\}\cap \{X_{T}\in B_{r}(y_{0})\}$ we
have
\begin{align*}
X_{T}-X_{T}^{\delta }= &\sum_{j=1}^{N}\int_{T-\delta }^{T}(\sigma
_{j}(t,X_{t})-\sigma _{j}(T-\delta ,X_{T-\delta }))1_{\{X_{T-\delta }\in
B_{2r}(y_{0})\}\cap \{X_{t}\in B_{2r}(y_{0})\}}dW_{t}^{j}+ \\
&+\int_{T-\delta }^{T}b(t,X_{t})dt.
\end{align*}%
If we are on the set $\{X_{T-\delta }\in B_{2r}(y_{0})\}\cap \{X_{t}\in
B_{2r}(y_{0})\}$ we may use (\ref{D3}) and we obtain%
\begin{align*}
\left\vert \sigma _{j}(t,X_{t})-\sigma _{j}(T-\delta ,X_{T-\delta
})\right\vert \leq &C\Big(1_{\{\left\vert X_{t}-X_{T-\delta }\right\vert
\geq \delta ^{-1/4}\}}+ \\
&+\left\vert \ln \left\vert X_{t}-X_{T-\delta }\right\vert \right\vert
^{-(1+h)}1_{\{\left\vert X_{t}-X_{T-\delta }\right\vert <\delta
^{-1/4}\}}+\left\vert \ln\frac1\delta\right\vert ^{-(2+h)}\Big) \\
\leq &C1_{\{\left\vert X_{t}-X_{T-\delta }\right\vert \geq \delta
^{-1/4}\}}+2C\left\vert \ln \frac{1}{\delta }\right\vert ^{-(2+h)}.
\end{align*}%
It follows that

\begin{eqnarray*}
&&\E(\left\vert \int_{T-\delta }^{T}(\sigma _{j}(t,X_{t})-\sigma
_{j}(T-\delta ,X_{T-\delta }))dW_{t}^{j}\right\vert ^{2}1_{\{\eta >T\}\cap
\{\tau >T\}\cap \{X_{T}\in B_{r}(y_{0})\}}) \\
&\leq &\int_{T-\delta }^{T}\E(\left\vert \sigma _{j}(t,X_{t})-\sigma
_{j}(T-\delta ,X_{T-\delta })\right\vert ^{2}1_{\{X_{T-\delta }\in
B_{2r}(y_{0})\}\cap \{X_{t}\in B_{2r}(y_{0})\}})dt \\
&\leq &C\delta \left\vert \ln \frac{1}{\delta }\right\vert
^{-2(2+h)}+\int_{T-\delta }^{T}\P(\left\vert X_{t}-X_{T-\delta }\right\vert
\geq \delta ^{-1/4})dt\leq C\delta \left\vert \ln \frac{1}{\delta }%
\right\vert ^{-2(2+h)}+C\delta ^{2}.
\end{eqnarray*}%
The drift term may be upper bounded directly:%
\begin{equation*}
\E(\left\vert \int_{T-\delta }^{T}b(t,X_{t})dt\right\vert ^{2})\leq
C\E(\left\vert \int_{T-\delta }^{T}(1+\left\vert X_{t}\right\vert
)dt\right\vert ^{2})\leq C\delta ^{2}.
\end{equation*}%
So we have proved that%
\begin{equation*}
\left\vert R_{\delta }\right\vert \leq \frac{C}{r}(\left\Vert \phi
\right\Vert _{\infty }+\left\Vert \nabla \phi \right\Vert _{\infty })\frac{%
\delta ^{1/2}}{\left\vert \ln \delta \right\vert ^{2+h}}
\end{equation*}%
and finally%
\begin{equation*}
\left\vert \int \phi d\nu -\int \phi d\nu _{\delta }\right\vert \leq \frac{C%
}{r}(\left\Vert \phi \right\Vert _{\infty }+\left\Vert \nabla \phi
\right\Vert _{\infty })\frac{\delta ^{1/2}}{\left\vert \ln \delta
\right\vert ^{2+h}}.
\end{equation*}%
$\square $

\smallskip

\textbf{Proof of Theorem \ref{9}.} \textbf{Step 1}. We define $\overline{\nu
}$ and $\overline{\nu }_{\delta }$ by%
\begin{equation*}
\int \phi d\overline{\nu }=\E(\phi \psi (X_{T\wedge \tau })1_{\{X_{T-\delta
}\in B_{2r}(y_{0})\}})\qquad and\qquad \int \phi d\overline{\nu }_{\delta
}=\E(\phi \psi (X_{T}^{\delta })1_{\{X_{T-\delta }\in B_{2r}(y_{0})\}}).
\end{equation*}%
One has
\begin{equation*}
\left\vert \int \phi d\overline{\nu }-\int \phi d\nu \right\vert \leq
\left\Vert \phi \psi \right\Vert _{\infty }\P(\{X_{T}\in B_{r}(y_{0})\}\cap
\{X_{T-\delta }\in B_{2r}^{c}(y_{0})\})\leq C\left\Vert \phi \psi
\right\Vert _{\infty }\delta
\end{equation*}%
so that $d_{1}(\nu ,\overline{\nu })\leq C\delta .$ In the same way $%
d_{1}(\nu _{\delta },\overline{\nu }_{\delta })\leq C\delta $ so that, using
the previous lemma
\begin{equation}
d_{1}(\nu ,\overline{\nu }_{\delta })\leq \frac{C}{r}\delta ^{1/2}(\ln
\delta )^{-(2+h)}  \label{diff1}
\end{equation}

\textbf{Step 2}. We have $\overline{\nu }_{\delta }(dy)=p_{\delta }(y)dy$
with%
\begin{equation*}
p_{\delta }(y)=\E(\psi (y)\gamma _{a_{\delta }(X_{T\wedge \tau
})}(y-X_{T\wedge \tau })1_{\{X_{T-\delta }\in B_{2r}(y_{0})\}})
\end{equation*}%
where $a_{\delta }(x)=\delta \sigma \sigma ^{\ast }(x)$ and
\begin{equation*}
\gamma _{a_{\delta }(x)}(y)=\frac{1}{\sqrt{\det a_{\delta }(x)}}\exp
(-\left\langle a_{\delta }^{-1}(x)y,y\right\rangle )
\end{equation*}%
is the Gaussian density of covariance matrix $a_{\delta }(x).$ We may find a
constant $\lambda _{\ast }$ such that the lower eigenvalue of $a_{\delta
}(x) $ is larger then $\delta \lambda _{\ast }$ for $x\in B_{2r}(y_{0}).$ It
is then easy to check that for every multi index $\alpha $ one has
\begin{equation*}
\sup_{x\in B_{2r}(y_{0})}\left\Vert \partial ^{\alpha }\gamma _{a_{\delta
}(x)}\right\Vert _{\infty }\leq C\delta ^{-\left\vert \alpha \right\vert /2}.
\end{equation*}%
We also have $\left\Vert \partial ^{\alpha }\psi \right\Vert _{\infty }\leq
Cr^{-\left\vert \alpha \right\vert }$ so it is easy to check that
\begin{equation*}
\left\Vert p_{\delta }\right\Vert _{2m+q,2m,p}\leq C(r,y_{0})\times \delta
^{-\frac{1}{2}(2m+q)}\qquad with\qquad C(r,y_{0})=\frac{C\left\vert
y_{0}\right\vert ^{2m}}{r^{2m+q}}.
\end{equation*}

\textbf{Step 3}. We will use the criterion (\ref{011}) so we focus on $%
e_{\log }.$ And we have for $p>1$%
\begin{equation*}
\left\Vert p_{\delta }\right\Vert _{2m+q,2m,e_{\log }}\leq \left\Vert
p_{\delta }\right\Vert _{2m+q,2m,p}\leq C(r,y_{0})\times \delta ^{-\frac{1}{2%
}(2m+q)}.
\end{equation*}%
We conclude that
\begin{equation*}
p_{\delta }\in M_{m,q,e_{\log }}(R_{\delta })\qquad with\qquad R_{\delta
}=C(r,y_{0})\times \delta ^{-\frac{1}{2}(2m+q)}.
\end{equation*}%
We have $\delta ^{1/2}=CR_{\delta }^{1/(2m+q)}$ so (\ref{diff1}) reads%
\begin{equation*}
d_{1}(\nu ,M_{m,q,e_{\log }}(R_{\delta }))\leq \frac{CR_{\delta }^{1/(2m+q)}%
}{(\ln R_{\delta })^{2+h}}.
\end{equation*}%
Now we look to (\ref{011}) with $k=1$ and $a>1$:%
\begin{equation*}
R_{\delta }^{\frac{q+1}{2m}}(\ln R_{\delta })^{a(1+\frac{q+1}{2m}%
)+1}d_{1}(\nu ,M_{m,q,e_{\log }}(R_{\delta }))\leq R_{\delta }^{\frac{q+1}{2m%
}-\frac{1}{2m+q}}(\ln R_{\delta })^{a(1+\frac{q+1}{2m})-1-h}.
\end{equation*}%
It is clear that if $q>0$ the above quantity blows up, but, if $q=0$ then $%
\frac{q+1}{2m}-\frac{1}{2m+q}=0$ and this term disappears. Moreover we take $%
m $ sufficiently large in order to have $\frac{1}{2m}<h$ and $a>1$
sufficiently small in order to have $a(1+\frac{1}{2m})\leq 1+h.$ And then
the above quantity is bounded. So $\nu $ is absolutely continuous and then $%
\mu $ is absolutely continuous on $B_{r}(y_{0}).$ $\square $

\subsection{A local H\"{o}rmander condition}

In this section we work under a local H\"{o}rmander condition. In contrast
with the situation from the previous section in this case we are no more
able to exhibit an explicit density (as the non degenerated Gaussian density
in the previous section) so we are obliged to use integration by parts. This
is the specific point in this example.

In this framework it is difficult to work with time dependent coefficients
so now on $\sigma _{j},b$ depend on $x$ only. Our equation is
\begin{equation*}
X_{t}=x+\sum_{j=1}^{N}\int_{0}^{t}\sigma
_{j}(X_{s})dW_{s}^{j}+\int_{0}^{t}b(X_{s})ds.
\end{equation*}%
We denote by $\sigma _{0}$ the drift coefficient when writing the equation
in Stratonovich form that is $\sigma _{0}=b-\frac{1}{2}\sum_{j=1}^{N}%
\sum_{k=1}^{d}\sigma _{j}^{k}\partial _{k}\sigma _{j}.$ We can do this only
if $\sigma _{j}$ are one time differentiable - in our case this will be
locally true, and this will be sufficient. For an open domain $D\subset
\R^{d} $ we denote by $C^{k}(D)$ the class of functions which are $k$ time
differentiable on $D.$ For $l\leq k$ we construct recursively%
\begin{equation*}
\mathcal{A}_{0}=\{\sigma _{j},j=1,...,N\},\quad \mathcal{A}_{l}=\{[\phi
,\sigma _{j}],\phi \in \mathcal{A}_{l-1},j=0,...,N\}
\end{equation*}%
and%
\begin{equation}
\Lambda _{k}(\sigma ,x)=\inf_{\left\vert \xi \right\vert =1}\sum_{\phi \in
\cup _{l=0}^{k}\mathcal{A}_{l}}\left\langle \phi (x),\xi \right\rangle ^{2}.
\label{D5}
\end{equation}%
So the condition $\Lambda _{0}(\sigma ,x)>0$ says that we have ellipticity
in the point $x$ and the condition $\Lambda _{k}(\sigma ,x)>0$ says that the
weak H\"{o}rmander condition of order $k$ holds in $x.$

\begin{theorem}\label{11}
Let $\Gamma \subset \R^{d}$ be an open domain and $\tau $ be the
exit time from $\Gamma$. Fix $y_{0}\in \Gamma $ and $r<\frac{%
1}{2}d(y_{0},\Gamma ^{c}).$ Suppose that the diffusion coefficients have linear
growth, $\sigma _{j},b\in C^{k}(B_{r}(y_{0})),j=0,1,...,N,$ for some $%
k\geq 2$ and moreover, $\Lambda _{k}(\sigma ,x)\geq \varepsilon
_{0}>0 $, for $x\in B_{r}(y_{0})$. Then for every $T>0$ the law $\mu $ of $%
X_{T\wedge \tau }$ is absolutely continuous with respect to the Lebesgue
measure with density $p_{\mu }\in C^{k-2}(B_{r/4}(y_{0})).$
\end{theorem}

\begin{remark}
In Proposition 23 from \cite{bib:[BCa]} we prove a similar result for a
diffusion process with coefficients which are globally in $C^{k}(\R^{d})$ and
for which we do not consider the stopping time $\tau .$ Here we generalize
this result for local regularity and with a stopping time.
\end{remark}

The proof is analogous to the proof of the previous theorem but now we will
use another approximation process. For $\delta >0$ we consider some
coefficients $\sigma _{j}^{\delta },b^{\delta }\in C^{k}(\R^{d}),j=1,...,N$
such that $\sigma _{j}^{\delta }(y)=\sigma _{j}(y),b^{\delta }(y)=b(y)$ for $%
y\in B_{r}(y_{0})$ and such that
\begin{equation}
\sum_{0\leq \left\vert \alpha \right\vert \leq k}\left\Vert \partial
^{\alpha }\sigma _{j}^{\delta }\right\Vert _{\infty }\leq \sum_{0\leq
\left\vert \alpha \right\vert \leq k}\left\Vert 1_{B_{r}(y_{0})}\partial
^{\alpha }\sigma _{j}\right\Vert _{\infty }  \label{A}
\end{equation}%
for $j=1,...,N$ and the same for $b^{\delta }.$ We also assume that%
\begin{equation}
\inf_{x\in \R^{N}}\Lambda _{k}(\sigma ^{\delta },x)\geq \inf_{x\in
B_{r}(y_{0})}\Lambda _{k}(\sigma ,x)=:\varepsilon _{\ast }>0.  \label{B}
\end{equation}%
Then we define $X_{t}^{\delta }=X_{t}$ for $t\leq T-\delta $ and
\begin{equation*}
X_{t}^{\delta }=X_{(T-\delta )\wedge \tau }+\sum_{j=1}^{N}\int_{(T-\delta
)\wedge \tau _{\Gamma }}^{t}\sigma _{j}^{\delta }(X_{s}^{\delta
})dW_{s}^{j}+\int_{(T-\delta )\wedge \tau _{\Gamma }}^{t}b(X_{s}^{\delta
})ds\quad \mbox{for } (T-\delta )\wedge \tau \leq t\leq T.
\end{equation*}%
We denote by $\mu $ the law of $X_{T\wedge \tau }$ and by $\mu _{\delta }$
the law of $X_{T}^{\delta }.$ Moreover we consider a truncation function $%
\psi \in C^{\infty }$ such that $1_{B_{r/2}(y_{0})}\leq \psi \leq
1_{B_{r}(y_{0})}$ and such that for every multi index $\alpha $ and every $%
p\geq 1$ there exists a constant $C_{\alpha ,p}$ such that
\begin{equation}
\left\vert \frac{\partial ^{\alpha }\psi }{\psi }(x)\right\vert ^{p}\times
\psi (x)\leq \frac{C_{\alpha ,p}}{r^{\left\vert \alpha \right\vert }}.
\label{Ineg}
\end{equation}%
Such a function is constructed for example in \cite{bib:[BCa]}, Section 2.7.
We define $\nu =\psi \mu $ and $\nu _{\delta }=\psi \mu
_{\delta }.$

\begin{lemma}
\label{12}Under the hypotheses of Theorem \ref{11}, for every $p>1$ one has
\begin{equation*}
d_{1}(\nu ,\nu _{\delta })\leq \frac{C_{p}}{r^{p}}\delta ^{p/2}.
\end{equation*}
\end{lemma}

\textbf{Proof}. We define $\eta =\inf \{t\geq T-\delta :\left\vert
X_{t}-X_{T-\delta }\right\vert >\frac{r}{4}\},\eta _{\delta }=\inf \{t\geq
T-\delta :\left\vert X_{t}^{\delta }-X_{T-\delta }^{\delta }\right\vert >%
\frac{r}{4}\}$ and we write%
\begin{align*}
&\left\vert \int \phi d\nu -\int \phi d\nu _{\delta }\right\vert \leq
I+2\left\Vert \phi \right\Vert _{\infty }(\P(\eta \leq T)+\P(\eta _{\delta
}\leq T))\qquad \text{with} \\
&I =\E(\left\vert \phi \psi (X_{T\wedge \tau })-\phi \psi (X_{T}^{\delta
})\right\vert 1_{\{\eta >T\}\cap \{\eta _{\delta }>T\}}).
\end{align*}%
Since the coefficients $\sigma _{j}$ have linear growth we have for every $p$%
\begin{equation*}
\P(\eta \leq T)\leq \P(\sup_{T-\delta \leq t\leq T}\left\vert
X_{t}(x)-X_{T-\delta }(x)\right\vert \geq \frac{r}{2})\leq \frac{C_{p}}{r^{p}%
}\delta ^{p/2}.
\end{equation*}%
And the same is true for $\P(\eta _{\delta }\leq T).$

We will now check that $I=0$. We write%
\begin{eqnarray*}
I &=&I^{\prime }+I^{\prime \prime }\qquad with \\
I^{\prime } &=&\E(\left\vert \phi \psi (X_{T\wedge \tau })-\phi \psi
(X_{T}^{\delta })\right\vert 1_{\{\eta >T\}\cap \{\eta _{\delta
}>T\}}1_{\{\tau <T}), \\
I^{\prime \prime } &=&\E(\left\vert \phi \psi (X_{T\wedge \tau })-\phi \psi
(X_{T}^{\delta })\right\vert 1_{\{\eta >T\}\cap \{\eta _{\delta
}>T\}}1_{\{\tau >T}).
\end{eqnarray*}%
We have $\psi (X_{\tau })=0.$ Suppose now that $\psi (X_{T}^{\delta })\neq 0$
and $\eta _{\delta }>T.$ Then $X_{(T-\delta )\wedge \tau }\in
B_{3r/4}(y_{0})\subset \Gamma $ so $\tau >T-\delta .$ But then, since $\eta
>T$ it follows that $\tau >T$ which is in contradiction with $\tau <T.$ So
we have proved that $I^{\prime }=0.$

We check now that $I^{\prime \prime }=0.$ We are on the set $\tau >T$ so $%
\tau $ disappears from the equation of $X_{t}^{\delta }.$ We are also on the
set $\{\eta >T\}.$\ If $\psi (X_{T}(x))\neq 0$\ it follows that $X_{T}\in
B_{r/4}(y_{0})$ and consequently for every $T-\delta \leq t\leq T$ we have $%
X_{t}\in B_{r/2}(y_{0}).$ On this set the coefficients $\sigma _{j}$ and $%
\sigma _{j}^{\delta }$ coincide so $X_{t}$ and $X_{t}^{\delta }$ solve the
same equation. Since the coefficients are Lipschitz continuous on $%
B_{r/2}(y_{0})$ this equation has the uniqueness property. We conclude that $%
X_{t}=X_{t}^{\delta }$ for $T-\delta \leq t\leq T.$ Suppose now that $\psi
(X_{T}^{\delta })\neq 0.$ Then, since we are on the set $\{\eta _{\delta
}>T\},$ the same reasoning gives that $X_{t}=X_{t}^{\delta }$ for $T-\delta
\leq t\leq T.$ So $I^{\prime \prime }=0$ and we obtain
\begin{equation*}
\left\vert \int \phi d\nu -\int \phi d\nu _{\delta }\right\vert \leq \frac{%
C_{p}}{r^{p}}\left\Vert \phi \right\Vert _{\infty }\delta ^{p/2}.
\end{equation*}%
$\square $

\smallskip

We now prove the following estimate for the Sobolev norms of $\nu
_{\delta }:$

\begin{lemma}
\label{13}Under the hypotheses of the Theorem \ref{11}, for every $q\leq k-1$
and $p>1$ there exist $C_{q,p}$ and $l_{q}$ such that
\begin{equation*}
\left\Vert 1\right\Vert _{W_{\nu _{\delta }}^{q,p}}\leq C_{q,p}\delta
^{-l_{q}}.
\end{equation*}
\end{lemma}

\textbf{Proof}. We consider the semigroup associated to the diffusion
process $X_{t}$ that is $P_{t}f(x)=\E(f(X_{t}(x)))$ where $X_{t}(x)$ is the
solution to our SDE starting from $X_{0}=x.$ And similarly $%
P_{t}^{\delta }f(x)=\E(f(X_{t}^{\delta }(x))).$ Then%
\begin{equation*}
\int fd\nu _{\delta }=\E(f(X_{T}^{\delta })\psi (X_{T}^{\delta }))=\int
P_{T-\delta }(x,dy)\int f(z)\psi (z)P_{\delta }^{\delta }(y,dz).
\end{equation*}%
A standard result concerning diffusion processes under H\"{o}rmander
conditions (see \cite{bib:[KS]}) says that for every multi index $%
\alpha $
\begin{equation*}
\E(\psi (X_{\delta }^{\delta }(x))\partial ^{\alpha }f(X_{\delta }^{\delta
}(x))=\E(f(X_{\delta }^{\delta })H_{\alpha }(X_{\delta }^{\delta }(x),\psi
(X_{\delta }^{\delta }(x))))
\end{equation*}%
where $H_{\alpha }(X_{\delta }^{\delta }(x),\psi (X_{\delta }^{\delta }(x)))$
is the weight which appears in the integration by parts formula using
Malliavin calculus. An inspection of the structure of this weight shows that
\begin{equation*}
H_{\alpha }(X_{\delta }^{\delta }(x),\psi (X_{\delta }^{\delta
}(x))=\sum_{\left\vert \beta \right\vert \leq \left\vert \alpha \right\vert
}\partial ^{\beta }\psi (X_{\delta }^{\delta }(x))\theta _{\beta }(x)
\end{equation*}%
where $\theta _{\beta }(x)$ is built using Malliavin derivatives of $%
X_{\delta }^{\delta }(x),$ the Ornstein Uhlembeck operator and the inverse
of the covariance matrix of $X_{\delta }^{\delta }(x).$ Moreover, under the
uniform H\"{o}rmander condition \ (\ref{B})\ and the boundedness condition (%
\ref{A})\ one can see that for every $p\geq 1$
\begin{equation*}
\E\left\vert \theta _{\beta }(x)\right\vert ^{p}\leq \frac{C}{\delta ^{l}}
\end{equation*}%
where $C$ depends on $p,$ on $\varepsilon _0$ and on the bounds in (%
\ref{A}) and $l\in \N$ is a power which depends on the order of the H\"{o}%
rmander condition. Moreover, taking conditional expectations we obtain
\begin{eqnarray*}
\E(f(X_{\delta }^{\delta }(x))H_{\alpha }(X_{\delta }^{\delta }(x),\psi
(X_{\delta }^{\delta }(x)))) &=&\E(f(X_{\delta }^{\delta }(x))\E(H_{\alpha
}(X_{\delta }^{\delta }(x),\psi (X_{\delta }^{\delta }(x))\mid X_{\delta
}^{\delta }(x))) \\
&=&\sum_{\left\vert \beta \right\vert \leq \left\vert \alpha \right\vert
}\E(f(X_{\delta }^{\delta }(x))\partial ^{\beta }\psi (X_{\delta }^{\delta
}(x))\E(\theta _{\beta }(x)\mid X_{\delta }^{\delta }(x))) \\
&=&\sum_{\left\vert \beta \right\vert \leq \left\vert \alpha \right\vert
}\E(f(X_{\delta }^{\delta }(x))\partial ^{\beta }\psi (X_{\delta }^{\delta
}(x))\widehat{\theta }_{\beta }(X_{\delta }^{\delta }(x)))
\end{eqnarray*}%
where $\widehat{\theta }_{\beta }(X_{\delta }^{\delta }(x))=\E(\theta _{\beta
}(x)\mid X_{\delta }^{\delta }(x)).$ Using H\"{o}lder inequality%
\begin{equation}
\E\left\vert \widehat{\theta }_{\beta }(x)\right\vert ^{p}\leq \E\left\vert
\theta _{\beta }(x)\right\vert ^{p}\leq \frac{C}{\delta ^{l}}.  \label{C}
\end{equation}%
We now come back and we write%
\begin{eqnarray*}
\int \partial ^{\alpha }fd\nu _{\delta } &=&\E(\partial ^{\alpha
}f(X_{T}^{\delta })\psi (X_{T}^{\delta }))=\int P_{T-\delta }(x,dy)\int
\partial ^{\alpha }f(z)\psi (z)P_{\delta }^{\delta }(y,dz) \\
&=&\int P_{T-\delta }(x,dy)\E(\partial ^{\alpha }f(X_{\delta }^{\delta
}(y))\psi (X_{\delta }^{\delta }(y))) \\
&=&\sum_{\left\vert \beta \right\vert \leq \left\vert \alpha \right\vert
}\int P_{T-\delta }(x,dy)\E(f(X_{\delta }^{\delta }(y))\partial ^{\beta }\psi
(X_{\delta }^{\delta }(y))\widehat{\theta }_{\beta }(X_{\delta }^{\delta
}(y)))) \\
&=&\sum_{\left\vert \beta \right\vert \leq \left\vert \alpha \right\vert
}\int P_{T-\delta }(x,dy)\E(f(X_{\delta }^{\delta }(y))\frac{\partial ^{\beta
}\psi (X_{\delta }^{\delta }(y))}{\psi (X_{\delta }^{\delta }(y))}\widehat{%
\theta }_{\beta }(X_{\delta }^{\delta }(y))\psi (X_{\delta }^{\delta }(y))))
\\
&=&\sum_{\left\vert \beta \right\vert \leq \left\vert \alpha \right\vert
}\int P_{T-\delta }(x,dy)\int f(z)\frac{\partial ^{\beta }\psi (z)}{\psi (z)}%
\widehat{\theta }_{\beta }(z)\psi (z)P_{\delta }^{\delta }(y,dz) \\
&=&\int f(z)\left( \sum_{\left\vert \beta \right\vert \leq \left\vert \alpha
\right\vert }\frac{\partial ^{\beta }\psi (z)}{\psi (z)}\widehat{\theta }%
_{\beta }(z)\right) d\nu _{\delta }(z).
\end{eqnarray*}%
This proves that%
\begin{equation*}
\partial _{\alpha }^{\nu _{\delta }}1(z)=\sum_{\left\vert \beta \right\vert
\leq \left\vert \alpha \right\vert }\frac{\partial ^{\beta }\psi (z)}{\psi
(z)}\widehat{\theta }_{\beta }(z).
\end{equation*}%
Now we compute%
\begin{eqnarray*}
\int \left\vert \partial _{\alpha }^{\nu _{\delta }}1(z)\right\vert ^{p}d\nu
_{\delta }(z) &\leq &C\sum_{\left\vert \beta \right\vert \leq \left\vert
\alpha \right\vert }\int \left\vert \frac{\partial ^{\beta }\psi (z)}{\psi
(z)}\widehat{\theta }_{\beta }(z)\right\vert ^{p}d\nu _{\delta }(z) \\
&=&C\sum_{\left\vert \beta \right\vert \leq \left\vert \alpha \right\vert
}\int \left\vert \frac{\partial ^{\beta }\psi (z)}{\psi (z)}\widehat{\theta }%
_{\beta }(z)\right\vert ^{p}\psi (z)d\mu _{\delta }(z) \\
&\leq &C\sum_{\left\vert \beta \right\vert \leq \left\vert \alpha
\right\vert }\int \left\vert \widehat{\theta }_{\beta }(z)\right\vert
^{p}d\mu _{\delta }(z),
\end{eqnarray*}%
last inequality being a consequence of the property (\ref{Ineg}) for $%
\psi .$ And by (\ref{C})
\begin{equation*}
\int \left\vert \widehat{\theta }_{\beta }(z)\right\vert ^{p}d\mu _{\delta
}(z)=\E\left\vert \widehat{\theta }_{\beta }(x)\right\vert ^{p}\leq \frac{C}{%
\delta ^{l}}.
\end{equation*}%
$\square $

\smallskip

\textbf{Proof of Theorem \ref{11}.} Since $\nu $ coincides with $\mu $ on $%
B_{r/2}(y_{0})$ we may look to the regularity of $\nu $ (instead of $\mu ).$
We recall now the hypothesis $\widetilde{H}_{q}(k,m,e)$ with $k=1$ (we work
with $d_{1}),m=1$ and $e(t)=t^{p}$ -in particular we have $\beta
_{e_{p}}(t)=t^{1/p_{\ast }}$ where $p_{\ast }$ is the conjugate of $p.$ Then, we consider hypothesis
\begin{align*}
\widetilde{H}_{q}(1,1,e_{p}):\quad &\overline{\lim }_{R\rightarrow \infty }%
\frac{L_{a}(R)^{1+\frac{1+q}{2}}L_{a}(R)^{\frac{d}{2p_{\ast }}}}{R}d_{1}(\nu
,\widetilde{M}_{1,q,e_{p}}(R))<\infty .
\end{align*}%
We also recall that $L_{a}(R)=R(\ln R)^{a}$ for $a>1$ so that%
\begin{equation*}
\frac{L_{a}(R)^{1+\frac{1+q}{2}}L_{a}(R)^{\frac{d}{2p_{\ast }}}}{R}=R^{\frac{%
1+q}{2}+\frac{d}{2p_{\ast }}}(\ln R)^{a(1+\frac{1+q}{2}+\frac{d}{2p_{\ast }}%
)}.
\end{equation*}%
We also recall that $\widetilde{M}_{1,q,e_{p}}(R)$ is the class of
probability measures $\rho $ such that $\left\Vert 1\right\Vert _{W_{\rho
}^{q,p}}\leq R.$ So in the previous lemmas we have proved that $\nu _{\delta
}\in \widetilde{M}_{1,q,e_{p}}(R_{\delta })$ with $R_{\delta }=\delta
^{-l_{q}}$ and $d_{1}(\nu ,\nu _{\delta })\leq C\delta ^{n}=R_{\delta
}^{-n/l_{q}}$ for some $l_{q}$ depending on $q$ and every $n\in \N.$ So we
have
\begin{equation*}
\frac{L_{a}(R)^{1+\frac{1+q}{2}}L_{a}(R)^{\frac{d}{2p_{\ast }}}}{R}d_{1}(\nu
,\widetilde{M}_{1,q,e_{p}}(R))\leq R^{\frac{1+q}{2}+\frac{d}{2p_{\ast }}%
}(\ln R)^{a(1+\frac{1+q}{2}+\frac{d}{2p_{\ast }})}\times R^{-n/l_{q}}
\end{equation*}%
and this inequality is true for every $n\in \N.$ So, by Theorem \ref{2F} we
have $\nu (dx)=p_{\nu }(x)dx$ with $p_{\nu }\in W^{k-1,p}$ for every $p\geq
1.$ Then using the Sobolev embedding theorem we obtain that $p_{\nu }\in
C^{k-2}.$ $\square $

\section{Stochastic heat equation}

\label{sect-heat}

In this section we investigate the regularity of the law of the solution to
the stochastic heat equation introduced by Walsh in \cite{bib:[W]}. Formally
this equation is%
\begin{equation}
\partial _{t}u(t,x)=\partial _{x}^{2}u(t,x)+\sigma (u(t,x))W(t,x)+b(u(t,x))
\label{S1}
\end{equation}%
where $W$ denotes a white noise on $\R_{+}\times [0,1].$ We consider
Neumann boundary conditions that is $\partial _{x}u(t,0)=\partial
_{x}u(t,1)=0 $ and the initial condition is $u(0,x)=u_{0}(x).$ The rigorous
formulation to this equation is given by the mild form constructed as
follows. Let $G_{t}(x,y)$ be the fundamental solution to the deterministic
heat equation $\partial _{t}v(t,x)=\partial _{x}^{2}v(t,x)$ with Neuman
boundary conditions. Then $u$ satisfies%
\begin{eqnarray}
u(t,x)
&=&\int_{0}^{1}G_{t}(x,y)u_{0}(y)dy+\int_{0}^{t}\int_{0}^{1}G_{t-s}(x,y)%
\sigma (u(s,y))dW(s,y)  \label{S2} \\
&&+\int_{0}^{t}\int_{0}^{1}G_{t-s}(x,y)b(u(s,y))dyds  \notag
\end{eqnarray}%
where $dW(s,y)$ is the It\^{o} integral introduced by Walsh. The function $%
G_{t}(x,y)$ is explicitly known (see \cite{bib:[W]} or \cite{bib:[BP]}) but
here we will use just few properties that we list below (see the appendix
in \cite{bib:[BP]} for the proof). More precisely, for $0<\varepsilon <t$ we
have%
\begin{equation}
\int_{t-\varepsilon }^{t}\int_{0}^{1}G_{t-s}^{2}(x,y)dyds\leq C\varepsilon
^{1/2}  \label{INEG1}
\end{equation}%
Moreover, for $0<x_{1}<...<x_{d}<1$ there exists a constant $C$ depending on
$\min_{i=1,d}(x_{i}-x_{i-1})$ such that%
\begin{equation}
C\varepsilon ^{1/2}\geq \inf_{\left\vert \xi \right\vert
=1}\int_{t-\varepsilon }^{t}\int_{0}^{1}\left( \sum_{i=1}^{d}\xi
_{i}G_{t-s}(x_{i},y)\right) ^{2}dyds\geq C^{-1}\varepsilon ^{1/2}.
\label{INEG2}
\end{equation}%
This is an easy consequence of the inequalities (A2) and (A3) from \cite%
{bib:[BP]}.

In \cite{bib:[PZ]} one gives sufficient conditions in order to obtain the
absolute continuity of the law of $u(t,x)$ for $(t,x)\in (0,\infty )\times
\lbrack 0,1]$ and in \cite{bib:[BP]}, under appropriate hypothesis, one
obtains a $C^{\infty }$ density for the law of the vector $%
(u(t,x_{1}),...,u(t,x_{d}))$ with $(t,x_{i})\in (0,\infty )\times \{\sigma
\neq 0\},i=1,...,d.$ The aim of this section is to obtain the same type of
results but under much weaker regularity hypothesis on the coefficients. One
may first discuss the absolute continuity of the law and further, under more
regularity hypothesis on the coefficients, one may discuss the regularity of
the density. Here, in order to avoid technicalities, we restrict ourself to
the absolute continuity property. We also assume global ellipticity that is
\begin{equation}\label{S3.1}
\sigma (x)\geq c_{\sigma }>0\qquad \text{for every } x\in [0,1].
\end{equation}%
A local ellipticity condition may also be used but again,\ this gives more
technical complications that we want to avoid. This is somehow a benchmark
for the efficiency of the method developed in the previous sections.

We assume the following regularity hypothesis: $\sigma ,b$ are measurable
and bounded functions and there exists $h>0$ such that

\begin{equation}
\left\vert \sigma (x)-\sigma (y)\right\vert \leq \left\vert \ln
\left\vert x-y\right\vert \right\vert ^{-(2+h)},\quad \text{for every }x,y\in [0,1].  \label{S3}
\end{equation}%
This hypothesis is not sufficient in order to ensure existence and
uniqueness for the solution to (\ref{S2}) (one needs $\sigma $ and $b$ to be
globally Lipschitz continuous in order to obtain it) - so in the following
we will just consider a random field $u(t,x),(t,x)\in (0,\infty )\times
\lbrack 0,1]$ which is adapted to the filtration generated by $W$ (see Walsh
\cite{bib:[W]} for precise definitions) and which solves (\ref{S2}).

\begin{proposition}
\label{SPDE1}
Suppose that (\ref{S3.1}) and (\ref{S3}) hold. Then for every $0<x_{1}<...<x_{d}<1$ and $T>0$, the law of the random vector $U=(u(T,x_{1}),...u(T,x_{d}))$ is absolutely continuous with respect to the Lebesgue measure.
\end{proposition}

\textbf{Proof}. Given $0<\varepsilon <T$ we decompose%
\begin{equation}
u(T,x)=u_{\varepsilon }(T,x)+I_{\varepsilon }(T,x)+J_{\varepsilon }(T,x)
\label{S4}
\end{equation}%
with%
\begin{align*}
u_{\varepsilon }(T,x)
=&\int_{0}^{1}G_{t}(x,y)u_{0}(y)dy+\int_{0}^{T}\int_{0}^{1}G_{T-s}(x,y)%
\sigma (u(s\wedge (T-\varepsilon ),y))dW(s,y) \\
&+\int_{0}^{T-\varepsilon }\int_{0}^{1}G_{T-s}(x,y)b(u(s,y))dyds,\\
I_{\varepsilon }(T,x)
=&
\int_{T-\varepsilon
}^{T}\int_{0}^{1}G_{T-s}(x,y)(\sigma (u(s,y))-\sigma (u(s\wedge
(T-\varepsilon ),y)))dW(s,y), \\
J_{\varepsilon }(T,x)
=&\int_{T-\varepsilon
}^{T}\int_{0}^{1}G_{T-s}(x,y)b(u(s,y))dyds.
\end{align*}
\textbf{Step 1}. We prove that
\begin{equation}
\E\left\vert I_{\varepsilon }(T,x)\right\vert ^{2}+\E\left\vert J_{\varepsilon
}(T,x)\right\vert ^{2}\leq C\left\vert \ln \varepsilon \right\vert
^{-2(2+h)}\varepsilon ^{1/2}.  \label{S5}
\end{equation}%
Let $\mu $ and $\mu_\varepsilon$ be the law of $%
U=(u(T,x_{1}),...,u(T,x_{d})) $ and $U_{\varepsilon }=(u_{\varepsilon
}(T,x_{1}),...,u_{\varepsilon }(T,x_{d}))$ respectively. Using the above
estimate one easily obtains
\begin{equation}
d_{1}(\mu ,\mu _{\varepsilon })\leq C\left\vert \ln \varepsilon \right\vert
^{-(2+h)}\varepsilon ^{1/4}.  \label{S5a}
\end{equation}%
Using the isometry property%
\begin{equation*}
\E\left\vert I_{\varepsilon }(T,x)\right\vert ^{2}=\int_{T-\varepsilon
}^{T}\int_{0}^{1}G_{T-s}^{2}(x,y)\E(\sigma (u(s,y)-\sigma (u(s\wedge
(T-\varepsilon ),y)))^{2})dyds.
\end{equation*}%
We consider the set $\Lambda _{\varepsilon ,\eta }(s,y)=\{\left\vert
u(s,y)-u(s\wedge (T-\varepsilon ),y)\right\vert \leq \eta \}$ and we split
the above term as $\E\left\vert I_{\varepsilon }(T,x)\right\vert
^{2}=A_{\varepsilon ,\eta }+B_{\varepsilon ,\eta }$ with

\begin{eqnarray*}
A_{\varepsilon } &=&\int_{T-\varepsilon
}^{T}\int_{0}^{1}G_{T-s}^{2}(x,y)\E(\sigma (u(s,y)-\sigma (u(s\wedge
(T-\varepsilon ),y)))^{2}1_{\Lambda _{\varepsilon ,\eta }(s,y)})dyds \\
B_{\varepsilon } &=&\int_{T-\varepsilon
}^{T}\int_{0}^{1}G_{T-s}^{2}(x,y)\E(\sigma (u(s,y)-\sigma (u(s\wedge
(T-\varepsilon ),y)))^{2}1_{\Lambda _{\varepsilon ,\eta }^{c}(s,y)})dyds.
\end{eqnarray*}%
Using (\ref{S3})%
\begin{equation*}
A_{\varepsilon }\leq C(\ln \eta )^{2(2+h)}\int_{T-\varepsilon
}^{T}\int_{0}^{1}G_{T-s}^{2}(x,y)dyds\leq C\left\vert \ln \eta \right\vert
^{-2(2+h)}\varepsilon ^{1/2}
\end{equation*}%
the last inequality being a consequence of (\ref{INEG1}). Moreover, coming
back to (\ref{S2}), we have
\begin{equation*}
\P(\Lambda _{\varepsilon ,\eta }^{c}(s,y))\leq \frac{1}{\eta ^{2}}\E\left\vert
u(s,y)-u(s\wedge (T-\varepsilon ),y)\right\vert ^{2}\leq \frac{C}{\eta ^{2}}%
\int_{T-\varepsilon }^{s}\int_{0}^{1}G_{s-r}^{2}(y,z)dzdr\leq \frac{%
C\varepsilon ^{1/2}}{\eta ^{2}}
\end{equation*}%
so that%
\begin{equation*}
B_{\varepsilon }\leq \frac{C\varepsilon ^{1/2}}{\eta ^{2}}%
\int_{T-\varepsilon }^{T}\int_{0}^{1}G_{T-s}^{2}(x,y)dyds\leq \frac{%
C\varepsilon }{\eta ^{2}}
\end{equation*}%
so that, taking $\eta =\varepsilon ^{1/16}$ we obtain%
\begin{equation*}
\E\left\vert I_{\varepsilon }(T,x)\right\vert ^{2}\leq C(\left\vert \ln
\varepsilon \right\vert ^{-2(2+h)}+\varepsilon ^{1/4})\varepsilon ^{1/2}\leq
C\left\vert \ln \varepsilon \right\vert ^{-2(2+h)}\varepsilon ^{1/2}.
\end{equation*}%
We estimate now%
\begin{equation*}
\left\vert J_{\varepsilon }(T,x)\right\vert \leq \left\Vert b\right\Vert
_{\infty }\int_{T-\varepsilon }^{T}\int_{0}^{1}G_{T-s}(x,y)dyds=\left\Vert
b\right\Vert _{\infty }\varepsilon
\end{equation*}%
so (\ref{S5}) is proved.

\textbf{Step 2}. Conditionally to $\mathcal{F}_{T-\varepsilon }$ the random
vector $U_{\varepsilon }=(u_{\varepsilon }(T,x_{1}),...,u_{\varepsilon
}(T,x_{d}))$ is Gaussian of covariance matrix%
\begin{equation*}
\Sigma ^{i,j}(U_{\varepsilon })=\int_{T-\varepsilon
}^{T}\int_{0}^{1}G_{T-s}(x_{i},y)G_{T-s}(x_{j},y)\sigma ^{2}(u(s\wedge
(T-\varepsilon ),y))dyds,\quad i,j=1,...,d.
\end{equation*}%
By (\ref{INEG2})
\begin{equation*}
C\sqrt{\varepsilon }\geq \Sigma (U_{\varepsilon })\geq \frac{1}{C}\sqrt{%
\varepsilon }
\end{equation*}%
where $C$ is a constant which depends on the upper bounds of $\sigma $ and
on $c_{\sigma }.$

We use now the criterion given in (\ref{O11}) with $k=1$ and $q=0.$ Let $%
p_{U_{\varepsilon }}$\ be the density of the law of $U_{\varepsilon }.$\
Conditionally to $\mathcal{F}_{T-\varepsilon }$ this is a Gaussian density
and direct estimates give (with the notation from Section \ref{sect-results}%
)
\begin{equation*}
\left\Vert p_{U_{\varepsilon }}\right\Vert _{2m,2m,e_{\log }}\leq
C\varepsilon ^{-m/2}.
\end{equation*}%
So if $\mu _{\varepsilon }$ is the law of $U_{\varepsilon }$ then $\mu
_{\varepsilon }\in M_{m,0,e_{\log }}(C\varepsilon ^{-m/2}).$ This is true
for every $m\in \N.$ Having in mind (\ref{S5a}) and taking $R_{\varepsilon
}=C\varepsilon ^{-m/2}$ the quantity in (\ref{O11}) is%
\begin{align*}
R_{\varepsilon }^{\frac{1}{2m}}\left\vert \ln R_{\varepsilon }\right\vert
^{a(1+\frac{1}{2m})+1}d_{1}(\mu ,\mu _{\varepsilon }) &=C\varepsilon
^{-1/4}\left\vert \ln \varepsilon \right\vert ^{a(1+\frac{1}{2m})+1}\times
\left\vert \ln \varepsilon \right\vert ^{-(2+h)}\varepsilon ^{1/4} \\
&=C\left\vert \ln \varepsilon \right\vert ^{a(1+\frac{1}{2m})-1-h}
\end{align*}%
where $h>0$ is fixed, being the one in (\ref{S3}), and $a>1$ is arbitrarily close
to one. Then taking $m$ sufficiently large we upper bound the above term by $%
C\left\vert \ln \varepsilon \right\vert ^{-h/2}$ and so we obtain
\begin{equation*}
\overline{\lim }_{\varepsilon \infty 0}R_{\varepsilon }^{\frac{1}{2m}%
}\left\vert \ln R_{\varepsilon }\right\vert ^{a(1+\frac{1}{2m})+1}d_{1}(\mu
,M_{m,0,e_{\log }}(R_{\varepsilon }))=0.
\end{equation*}%
Using now the result given in Example 2 we conclude that $\mu $ is
absolutely continuous. $\square $

\section{Appendix}

\subsection{Super kernels}

We consider a super kernel $\phi :\R^{d}\rightarrow \R_{+}$ that is a function
which has the support included in $B_{1}(0),$ $\int \phi (x)dx=1$ and such
that for every non null multi index $\alpha =(\alpha _{1},...,\alpha
_{d})\in \N^{d}$ one has
\begin{equation}
\int y^{\alpha }\phi (y)dy=0\qquad y^{\alpha }=\prod_{i=1}^{d}y_{i}^{\alpha
_{i}}.  \label{kk1}
\end{equation}%
See (\cite{[KK]}) Section 3, Remark 1 for the construction of a
superkernel.\ The corresponding $\phi _{\delta }$, $\delta \in (0,1)$, is
defined by
\begin{equation*}
\phi _{\delta }(y)=\frac{1}{\delta ^{d}}\phi (\frac{y}{\delta }).
\end{equation*}%
For a function $f$ we denote $f_{\delta }=f\ast \phi _{\delta }.$ We will
work with the norms $\left\Vert f\right\Vert _{k,\infty }$ and $\left\Vert
f\right\Vert _{q,l,(e)}$ defined in (\ref{O4}) and in (\ref{O4a}). And we
have

\begin{lemma}
\label{kernel copy(1)} i) Let $k,q\in \N,l>d$ and $e\in \mathcal{E}$. There
exists a universal constant $C$ such that for every $f\in W^{q,l,e}$ one has%
\begin{equation}
\left\Vert f-f_{\delta }\right\Vert _{W_{\ast }^{k,\infty }}\leq C\left\Vert
f\right\Vert _{q,l,(e)}\delta ^{q+k}.  \label{kk2}
\end{equation}%
ii) Let $k,q\in \N$ and let $e\in \mathcal{E}$. There exists a universal
constant $C$ such that for every $f\in W^{q,0,e}$ one has%
\begin{equation}
\left\Vert f-f_{\delta }\right\Vert _{W_{\ast }^{k,e_{\ast }}}\leq
C\left\Vert f\right\Vert _{q,0,(e)}\delta ^{q+k}.  \label{kk4}
\end{equation}%
iii) Let $l>d,n,q\in \N$, with $n\geq q$, and let $e\in \mathcal{E}$. There
exists a universal constant $C$ such that%
\begin{equation}
\left\Vert f_{\delta }\right\Vert _{n,l,(e)}\leq C\left\Vert f\right\Vert
_{q,l,(e)}\delta ^{-(n-q)}.  \label{kk3}
\end{equation}
\end{lemma}

\textbf{Proof}. i) We may suppose without loss of generality that $f\in
C_{b}^{\infty }.$ Using Taylor expansion of order $q+k$
\begin{eqnarray*}
f(\emph{x)}-f_{\delta }(x) &=&\int (f(\emph{x)}-f(y))\phi _{\delta }(x-y)dy
\\
&=&\int I(x,y)\phi _{\delta }(x-y)dy+\int R(x,y)\phi _{\delta }(x-y)dy
\end{eqnarray*}%
with
\begin{eqnarray*}
I(x,y) &=&\sum_{i=1}^{q+k-1}\frac{1}{i!}\sum_{\left\vert \alpha \right\vert
=i}\partial ^{\alpha }f(x)(x-y)^{\alpha }, \\
R(x,y) &=&\frac{1}{(q+k)!}\sum_{\left\vert \alpha \right\vert
=q+k}\int_{0}^{1}\partial ^{\alpha }f(x+\lambda (y-x))(x-y)^{\alpha
}d\lambda .
\end{eqnarray*}%
Using (\ref{kk1}) we obtain $\int I(x,y)\phi _{\delta }(x-y)dy=0$ and by a
change of variable we get
\begin{equation*}
\int R(x,y)\phi _{\delta }(x-y)dy=\frac{1}{(q+k)!}\sum_{\left\vert \alpha
\right\vert =q+k}\int_{0}^{1}\int dz\phi _{\delta }(z)\partial ^{\alpha
}f(x+\lambda z)z^{\alpha }d\lambda .
\end{equation*}%
We consider now $g\in W^{k,\infty }$ and we write%
\begin{equation*}
\int (f(\emph{x)}-f_{\delta }(x))g(x)dx=\frac{1}{(q+k)!}\sum_{\left\vert
\alpha \right\vert =q+k}\int_{0}^{1}d\lambda \int dz\phi _{\delta
}(z)z^{\alpha }\int \partial ^{\alpha }f(x+\lambda z)g(x)dx.
\end{equation*}%
Let us denote $f_{a}(x)=f(x+a).$ We have $(\partial ^{\alpha
}f)(x+a)=(\partial ^{\alpha }f_{a})(x).$ Let $\alpha $ with $\left\vert
\alpha \right\vert =\sum_{i=1}^{d}\alpha _{i}=q+k.$ We split $\alpha $ into
two multi indexes $\beta $ and $\gamma $ such that $\left\vert \beta
\right\vert =k,\left\vert \gamma \right\vert =q$ and $\partial ^{\beta
}\partial ^{\gamma }=\partial ^{\alpha }$ (this may be done in several ways
but any one of them is good for us). Then using integration by parts%
\begin{eqnarray*}
\left\vert \int \partial ^{\alpha }f(x+\lambda z)g(x)dx\right\vert
&=&\left\vert \int \partial ^{\beta }\partial ^{\gamma }f_{\lambda
z}(x)g(x)dx\right\vert \\
&\leq &\int \left\vert \partial ^{\gamma }f_{\lambda z}(x)\right\vert
\left\vert \partial ^{\beta }g(x)\right\vert dx\leq \left\Vert g\right\Vert
_{k,\infty }\int \left\vert \partial ^{\gamma }f_{\lambda z}(x)\right\vert dx
\\
&=&\left\Vert g\right\Vert _{k,\infty }\int \left\vert \partial ^{\gamma
}f(x)\right\vert dx.
\end{eqnarray*}%
We write $\partial ^{\gamma }f(x)=u_{l}(x)v_{\gamma }(x)$ with $%
u_{l}(x)=(1+\left\vert x\right\vert ^{2})^{l/2}$ and $v_{\gamma
}(x)=(1+\left\vert x\right\vert ^{2})^{-l/2}\partial ^{\alpha }f(x).$ Using H%
\"{o}lder inequality%
\begin{equation*}
\int \left\vert \partial ^{\gamma }f(x)\right\vert dx\leq C\left\Vert
u_{l}\right\Vert _{(e_{\ast })}\left\Vert v_{\gamma }\right\Vert _{(e)}\leq
C\left\Vert u_{l}\right\Vert _{(e_{\ast })}\left\Vert f\right\Vert
_{q,l,(e)}.
\end{equation*}%
By Remark \ref{U} $\left\Vert u_{l}\right\Vert _{(e_{\ast })}<\infty .$ So
we obtain
\begin{eqnarray*}
\left\vert \int_{0}^{1}\int dz\phi _{\delta }(z)z^{\alpha }\int \partial
^{\alpha }f(x+\lambda z)g(x)dxd\lambda \right\vert &\leq &C\left\Vert
f\right\Vert _{q,l,(e)}\left\Vert g\right\Vert _{k,\infty }\int \phi
_{\delta }(z)\left\vert z\right\vert ^{k+q}dz \\
&\leq &C\left\Vert f\right\Vert _{q,l,(e)}\left\Vert g\right\Vert _{k,\infty
}\delta ^{k+q}.
\end{eqnarray*}

ii) The proof is exactly the same but one uses directly\ H\"{o}lder's
inequality
\begin{equation*}
\int \left\vert \partial ^{\gamma }f_{\lambda z}(x)\right\vert \left\vert
\partial ^{\beta }g(x)\right\vert dx\leq 2\left\Vert \partial ^{\gamma
}f_{\lambda z}\right\Vert _{(e)}\left\Vert \partial ^{\beta }g\right\Vert
_{(e_{\ast })}.
\end{equation*}%
And we have%
\begin{eqnarray*}
\left\Vert \partial ^{\gamma }f_{\lambda z}\right\Vert _{(e)} &=&\inf
\{c:\int e(\frac{1}{c}\partial ^{\gamma }f(x+\lambda z))dx\leq 1\} \\
&=&\inf \{c:\int e(\frac{1}{c}\partial ^{\gamma }f(x))dx\leq 1\}=\left\Vert
\partial ^{\gamma }f\right\Vert _{(e)}.
\end{eqnarray*}

iii) Let $\alpha $ be a multi index with $\left\vert \alpha \right\vert =n$
and let $\beta ,\gamma $ be a splitting of $\alpha $ with $\left\vert \beta
\right\vert =q$ and $\left\vert \gamma \right\vert =n-q.$ We have
\begin{align*}
u(x)& :=(1+\left\vert x\right\vert )^{l}\left\vert \partial ^{\alpha
}f_{\delta }(x)\right\vert =(1+\left\vert x\right\vert )^{l}\left\vert
\partial ^{\beta }f\ast \partial ^{\gamma }\phi _{\delta }(x)\right\vert \\
& \leq (1+\left\vert x\right\vert )^{l}\left\vert \partial ^{\beta
}f\right\vert \ast \left\vert \partial ^{\gamma }\phi _{\delta }\right\vert
(x)\leq v\ast \left\vert \partial ^{\gamma }\phi _{\delta }\right\vert (x)
\end{align*}%
with $v(x)=2^{l}(1+\left\vert x\right\vert )^{l}\left\vert \partial ^{\beta
}f\right\vert (x).$ The last inequality is due to the fact that if $%
\left\vert \partial ^{\gamma }\phi _{\delta }\right\vert (x-y)\neq 0$ then$%
(1+\left\vert y\right\vert )^{l}\leq 2(1+\left\vert x\right\vert )^{l}.$
Using (\ref{Oo2})\ we obtain
\begin{equation*}
\left\Vert u\right\Vert _{(e)}\leq \left\Vert v\ast \left\vert \partial
^{\gamma }\phi _{\delta }\right\vert \right\Vert _{(e)}\leq \left\Vert
\partial ^{\gamma }\phi _{\delta }\right\Vert _{1}\left\Vert v\right\Vert
_{(e)}\leq \frac{C}{\delta ^{n-q}}\left\Vert v\right\Vert _{(e)}\leq \frac{C%
}{\delta ^{n-q}}\left\Vert f_{\beta ,l}\right\Vert _{(e)}.
\end{equation*}%
$\square $

\begin{proposition}
\label{ante-rec} Let $r,k,n\in \N$, with $n>r$ and $l>d$, and $e\in \mathcal{E%
}$. We take $\theta =\frac{k+r}{n-r}$ and $\gamma <\frac{\theta }{1+\theta }=%
\frac{k+r}{k+n}$. Then (with the notation from Section 3)%
\begin{equation}
W^{r,l,e}\subset B_{\theta ,0}(W^{n,l,e},W_{\ast }^{k,\infty })\subset
(W^{n,l,e},W_{\ast }^{k,\infty })_{\gamma }.  \label{K}
\end{equation}
\end{proposition}

\textbf{Proof}. We denote $X=W^{n,l,e},Y=W_{\ast }^{k,\infty }.$ Let $f\in
W^{r,l,e}$ and $\delta \in (0,1).$ From (\ref{kk3}) we have%
\begin{equation*}
\left\Vert f_{\delta }\right\Vert _{X}=\left\Vert f_{\delta }\right\Vert
_{n,l,(e)}\leq C\left\Vert f\right\Vert _{r,l,(e)}\delta ^{-(n-r)}=:R.
\end{equation*}%
So $f_{\delta }\in B_{X}(R).$ On the other hand, (\ref{kk2}) gives
\begin{equation*}
\left\Vert f-f_{\delta }\right\Vert _{Y}=\left\Vert f-f_{\delta }\right\Vert
_{W_{\ast }^{k,\infty }}\leq C\left\Vert f\right\Vert _{r,l,(e)}\delta
^{k+r}=C\left\Vert f\right\Vert _{r,l,(e)}^{\rho }R^{-\frac{k+r}{n-r}},
\end{equation*}%
so that
\begin{equation*}
R^{\frac{k+r}{n-r}}d_{Y}(f,B_{X}(R))\leq C\left\Vert f\right\Vert
_{r,l,(e)}^{\rho }.
\end{equation*}%
This means that $f\in B_{\theta ,0}(X,Y).$

We prove now the second inclusion. We have $B_{\theta ,0}(X,Y)\subset
B_{\alpha ,\beta }(X,Y)$ for every $\alpha <\theta $ and every $\beta \geq
0. $ We recall that by Proposition \ref{Spaces} we have $B_{\alpha ,\beta
}(X,Y)\subset (X,Y)_{\gamma }$ with $\alpha =\frac{\gamma }{1-\gamma }$ and $%
\beta =\frac{2}{1-\gamma }.$ By our hypothesis $\alpha =\frac{\gamma }{%
1-\gamma }<\theta $ so we obtain $B_{\theta ,0}(X,Y)\subset B_{\alpha ,\beta
}(X,Y)\subset (X,Y)_{\gamma }.$ $\square $

\smallskip

Recall now Hypothesis $H_{q}(k,m,e)$ and the set $B_{q}(k,m,e)$ defined in
Section \ref{sect-results}. To shorten notations, we put $%
X=W^{2m+q,2m,e},Y=W_{\ast }^{k,\infty }.$\ Then $M_{m,q,e}(R)=B_{X}(R)$ and $%
d_{k}(\mu ,M_{m,q,e}(R))=d_{Y}(\mu ,B_{X}(R)).$ For $a>1$ we have denoted $%
L_{a}(R)=R(\ln R)^{a}$ and Hypothesis $H_{q}(k,m,e)$ reads: there exists $%
a>1 $ such that%
\begin{equation}
\overline{\lim }_{R\rightarrow \infty }\frac{L_{a}(R)^{1+\frac{k+q}{2m}}}{R}%
\beta_e (L_{a}(R)^{\frac{d}{2m}})d_{Y}(u,B_{X}(R))<\infty .  \label{a2}
\end{equation}
And $B_{q}(k,m,e)$ is the set of measures $\mu $ such that $H_{q}(k,m,e)$
holds for $\mu .$

\begin{corollary}
\label{Reciproc}Let $k,q,m\in \N$ with $m>d/2$, and $e\in \mathcal{E}_{\alpha
,\beta }$ (see (\ref{RR})). Suppose that $0\leq \alpha <\frac{2m+k+q}{d(2m-1)%
}$ and $\beta \geq 0$. Then%
\begin{equation*}
W^{q+1,2m,e}\subset B_{q}(k,m,e).
\end{equation*}
\end{corollary}

\textbf{Proof}. We take $r=q+1,n=2m+q,l=2m$ in Proposition \ref{ante-rec}.
Then $\theta =\frac{k+q+1}{2m-1}.$ So we know that for $u\in W^{q+1,2m,e}$
we have%
\begin{equation}
\overline{\lim }_{R\rightarrow \infty }R^{\frac{k+q+1}{2m-1}%
}d_{Y}(u,B_{X}(R))<\infty .  \label{a1}
\end{equation}%
Moreover, by (\ref{RR}), for sufficiently large $R$
\begin{equation*}
\frac{L_{a}(R)^{1+\frac{k+q}{2m}}}{R}\beta _{e}(L_{a}(R)^{\frac{d}{2m}})\leq
CR^{\frac{k+q}{2m}+\frac{d\alpha }{2m}}(\ln R)^{\rho }
\end{equation*}%
with $\rho =a(1+\frac{k+q}{2m}+\frac{d}{2m})+\beta .$ Our hypothesis on $%
\alpha $ ensures that $\frac{k+q+1}{2m-1}>\frac{k+q}{2m}+\frac{d\alpha }{2m}$
and so (\ref{a1}) implies (\ref{a2}). $\square $

\subsection{Norms}

In this section we use the notation introduced in Section \ref{sect-interp}
which we recall here. We consider two normed spaces $(X,\left\Vert \circ
\right\Vert _{X})$ and $(Y,\left\Vert \circ \right\Vert _{Y})$ such that $%
X\subset Y.$ We also consider some $a\geq 0,m\in \N_{\ast }$ and $\theta >0.$
For $y\in Y$ and for a sequence $x_{n}\in X,n\in \N$ we define
\begin{equation*}
\pi _{\theta ,m,a}(y,(x_{n})_{n})=\sum_{n=1}^{\infty }2^{n\theta
}n^{a}\left\Vert y-x_{n}\right\Vert _{Y}+\frac{1}{2^{2nm}}\left\Vert
x_{n}\right\Vert _{X}.
\end{equation*}%
Moreover we define $\rho _{\theta ,m,a}^{X,Y}(y)=\inf \pi _{\theta ,m,a}(y)$
with the infimum taken over all the sequences $x_{n}\in X,n\in \N.$ Finally
we denote
\begin{equation*}
S_{\theta ,m,a}(X,Y)=\{y\in Y:\rho _{\theta ,m,a}^{X,Y}(y)<\infty \}.
\end{equation*}%
Moreover we denote $K(y,t)=\inf \{\left\Vert y-x\right\Vert _{Y}+t\left\Vert
x\right\Vert _{X}\}$ with the infimum taken over all $x\in X$ and we define%
\begin{equation*}
|y|_{\gamma ,b}=\int_{0}^{1}\frac{\left\vert \ln t\right\vert ^{b}}{%
t^{\gamma }}K(y,t)\frac{dt}{t}.
\end{equation*}%
We denote
\begin{equation*}
K_{\gamma ,b}(X,Y)=\{y\in Y:|y|_{\gamma ,b}<\infty \}.
\end{equation*}

\begin{proposition}
\label{NORM}We have
\begin{equation*}
S_{\theta ,m,a}(X,Y)=K_{\gamma ,b}(X,Y)\qquad with\qquad \gamma =\frac{%
\theta }{2m+\theta },b=\frac{2ma}{2m+\theta }.
\end{equation*}%
and there exists a universal constant $C$ (which may be computed explicitly)
such that
\begin{equation*}
\frac{1}{C}\rho^{X,Y} _{\theta ,m,a}(y)\leq | y| _{\gamma ,b}\leq
C(\left\Vert y\right\Vert _{Y}+\rho^{X,Y} _{\theta ,m,a}(y)).
\end{equation*}
\end{proposition}

\textbf{Proof}. \textbf{Step 1}. We write%
\begin{equation*}
\pi _{\theta ,m,a}(y,(x_{n})_{n})=\sum_{n}2^{n\theta }n^{a}(\left\Vert
y-x_{n}\right\Vert _{Y}+\frac{1}{n^{a}2^{n(2m+\theta )}}\left\Vert
x_{n}\right\Vert _{X})
\end{equation*}%
and we define%
\begin{equation*}
t_{n}=\frac{1}{n^{a}2^{n(2m+\theta )}}.
\end{equation*}%
We have%
\begin{equation*}
t_{n}-t_{n+1}=t_{n}\alpha _{n}\qquad with\qquad \alpha _{n}=1-\frac{1}{%
2^{2m+\theta }}\times \left( \frac{n}{n+1}\right) ^{a}.
\end{equation*}%
And%
\begin{equation*}
\alpha ^{\ast }:=1-\frac{1}{2^{2m+\theta +a}}\geq \alpha _{n}\geq 1-\frac{1}{%
2^{2m+\theta }}=:\alpha _{\ast }.
\end{equation*}%
Then we write%
\begin{equation*}
\pi _{\theta ,m,a}(y,(x_{n})_{n})=\sum_{n}2^{n\theta }n^{a}\times \frac{%
\left\Vert y-x_{n}\right\Vert _{Y}+t_{n}\left\Vert x_{n}\right\Vert _{X}}{%
t_{n}\alpha _{n}}\times (t_{n}-t_{n+1})
\end{equation*}%
so that%
\begin{equation*}
\frac{1}{\alpha ^{\ast }}\pi _{\theta ,m,a}^{\prime }(y,(x_{n})_{n})\leq \pi
_{\theta ,m,a}(y,(x_{n})_{n})\leq \frac{1}{\alpha _{\ast }}\pi _{\theta
,m,a}^{\prime }(y,(x_{n})_{n})
\end{equation*}%
with%
\begin{equation*}
\pi _{\theta ,m,a}^{\prime }(y,(x_{n})_{n})=\sum_{n}2^{n\theta }n^{a}\times
\frac{\left\Vert y-x_{n}\right\Vert _{Y}+t_{n}\left\Vert x_{n}\right\Vert
_{X}}{t_{n}}\times (t_{n}-t_{n+1}).
\end{equation*}

\textbf{Step 2}. We have
\begin{equation*}
\frac{\left\vert \ln t_{n}\right\vert ^{b}}{t_{n}^{\theta /(2m+\theta )}}%
=\left( a\ln n+n(2m+\theta )\ln 2\right) ^{b}n^{a\theta /(2m+\theta
)}2^{n\theta }.
\end{equation*}%
Since $b+a\theta /(2m+\theta )=a$ we may find $C$ such that%
\begin{equation*}
C\times n^{a}2^{n\theta }\geq \frac{\left\vert \ln t_{n}\right\vert ^{b}}{%
t_{n}^{\theta /(2m+\theta )}}\geq \frac{1}{C}\times n^{a}2^{n\theta }.
\end{equation*}%
Notice that the functions $t\rightarrow t^{-1}K(y,t)$ \ and $t\rightarrow
t^{-^{1}y/(2m+\theta )}\left\vert \ln t\right\vert ^{b}$\ are decreasing so
we obtain
\begin{equation*}
\frac{\left\vert \ln t_{n+1}\right\vert ^{b}}{t_{n+1}^{\theta /(2m+\theta )}}%
\frac{K(y,t_{n+1})}{t_{n+1}}(t_{n}-t_{n+1})\geq \int_{t_{n+1}}^{t_{n}}\frac{%
\left\vert \ln t\right\vert ^{b}}{t^{\theta /(2m+\theta )}}K(y,t)\frac{dt}{t}%
\geq \frac{\left\vert \ln t_{n}\right\vert ^{b}}{t_{n}^{\theta /(2m+\theta )}%
}\frac{K(y,t_{n})}{t_{n}}(t_{n}-t_{n+1}).
\end{equation*}%
By the very definition of $K(y,t_{n})$ we may find $x_{n}\in X$ such that $%
K(y,t_{n})\geq \frac{1}{2}(\left\Vert y-x_{n}\right\Vert
_{Y}+t_{n}\left\Vert x_{n}\right\Vert _{X}).$ It follows that
\begin{eqnarray*}
\int_{0}^{1}\frac{\left\vert \ln t\right\vert ^{b}}{t^{\theta /(2m+\theta )}}%
K(y,t)\frac{dt}{t} &\geq &\sum_{n}\int_{t_{n+1}}^{t_{n}}\frac{\left\vert \ln
t\right\vert ^{b}}{t^{\theta /(2m+\theta )}}K(y,t)\frac{dt}{t} \\
&\geq &\frac{1}{2}\sum_{n}\frac{\left\vert \ln t_{n}\right\vert ^{b}}{%
t_{n}^{\theta /(2m+\theta )}}\frac{\left\Vert y-x_{n}\right\Vert
_{Y}+t_{n}\left\Vert x_{n}\right\Vert _{X}}{t_{n}}(t_{n}-t_{n+1}) \\
&\geq &\frac{1}{2C}\pi _{\theta ,m,a}^{\prime }(y,(x_{n})_{n}).
\end{eqnarray*}%
So we have proved that for $\gamma =\theta /(2m+\theta )$ and $%
b=2ma/(2m+\theta )$ one has
\begin{equation*}
|y|_{\gamma ,b}\geq \frac{\alpha _{\ast }}{2C}\pi _{\theta
,m,a}(y,(x_{n})_{n})\geq \frac{\alpha _{\ast }}{2C}\rho _{\theta
,m,a}^{X,Y}(y).
\end{equation*}

We write now
\begin{eqnarray*}
t_{n}-t_{n+1} &=&\alpha _{n}t_{n}=\frac{\alpha _{n}}{\alpha _{n+1}}\times
\frac{t_{n}}{t_{n+1}}\times \alpha _{n+1}t_{n+1}=\frac{\alpha _{n}}{\alpha
_{n+1}}\times \frac{t_{n}}{t_{n+1}}\times (t_{n+1}-t_{n+2}) \\
&\leq &2^{2m+\theta +a+1}(t_{n+1}-t_{n+2}).
\end{eqnarray*}%
Then%
\begin{equation*}
2^{2m+\theta +a+1}\times \frac{\left\vert \ln t_{n+1}\right\vert ^{b}}{%
t_{n+1}^{\theta /(2+\theta )}}\frac{K(y,t_{n+1})}{t_{n+1}}%
(t_{n+1}-t_{n+2})\geq \int_{t_{n+1}}^{t_{n}}\frac{\left\vert \ln
t\right\vert ^{b}}{t^{\theta /(2+\theta )}}K(y,t)\frac{dt}{t}.
\end{equation*}%
For every $x_{n+1}\in X$we have $K(y,t_{n+1})\leq \left\Vert
y-x_{n+1}\right\Vert _{Y}+t_{n+1}\left\Vert x_{n+1}\right\Vert _{X}$ so that
for every sequence $x_{n}\in X$ we obtain
\begin{align*}
&2^{2m+\theta +a+1}\times \sum_{n}\frac{\left\vert \ln t_{n+1}\right\vert ^{b}%
}{t_{n+1}^{\theta /(2+\theta )}}\frac{\left\Vert y-x_{n+1}\right\Vert
_{Y}+t_{n+1}\left\Vert x_{n+1}\right\Vert _{X}}{t_{n+1}}(t_{n+1}-t_{n+2})%
\geq\\
&\qquad\qquad\geq  \sum_{n}\int_{t_{n+1}}^{t_{n}}\frac{\left\vert \ln t\right\vert ^{b}}{%
t^{\theta /(2+\theta )}}K(y,t)\frac{dt}{t}.
\end{align*}%
This means that
\begin{equation*}
2^{2m+\theta +a+1}\alpha ^{\ast }\pi _{\theta ,m,a}(y,(x_{n})_{n})\geq
2^{2m+\theta +a+1}\pi _{\theta ,m,a}^{\prime }(y,(x_{n})_{n})\geq
\int_{0}^{t_{1}}\frac{\left\vert \ln t\right\vert ^{b}}{t^{\theta /(2+\theta
)}}K(y,t)\frac{dt}{t}.
\end{equation*}%
Since this inequality holds for every sequence it holds for the infimum
also. So we obtain%
\begin{equation*}
2^{2m+\theta +a+1}\alpha ^{\ast }\rho _{\theta ,m,a}^{X,Y}(y)\geq
\int_{0}^{t_{1}}\frac{\left\vert \ln t\right\vert ^{b}}{t^{\theta /(2+\theta
)}}K(y,t)\frac{dt}{t},
\end{equation*}%
and the statement follows. $\square $

\smallskip

We define now%
\begin{equation*}
B_{\alpha ,\beta }(X,Y)=\{y\in Y:\overline{\lim }_{R\rightarrow \infty
}R^{\alpha }(\ln R)^{\beta }d_{Y}(y,B_{X}(R))<\infty \}
\end{equation*}

\begin{proposition}
\label{balance}%
\begin{equation*}
B_{\alpha ,\beta }(X,Y)\subset S_{\theta ,m,a}(X,Y)\qquad with\qquad \alpha =%
\frac{\theta }{2m},\beta =2+a+\frac{\theta }{m}.
\end{equation*}
\end{proposition}

\textbf{Proof}. If $y\in B_{\alpha ,\beta }(X,Y)$ one may find $R_{\ast
},C_{\ast }$ such that for every $R\geq R_{\ast }$ there exists $x_{R}\in X$
such that%
\begin{equation*}
\left\Vert x_{R}\right\Vert _{X}\leq R\qquad and\qquad \left\Vert
y-x_{R}\right\Vert _{X}\leq \frac{C_{\ast }}{R^{\alpha }(\ln R)^{\beta }}.
\end{equation*}%
We take $R_{n}=n^{-2}2^{2nm}$ and $n_{\ast }$ such that $R_{n_{\ast }}\geq
R_{\ast }.$ Then%
\begin{equation*}
\sum_{n\geq n_{\ast }}\frac{1}{2^{2nm}}\left\Vert x_{Rn}\right\Vert _{X}\leq
\sum_{n\geq n_{\ast }}\frac{1}{n^{2}}<\infty .
\end{equation*}%
Moreover, since $2m\alpha =\theta $ and $\beta -a-2\alpha =2$ we have%
\begin{equation*}
n^{a}2^{n\theta }\times \frac{1}{R_{n}^{\alpha }(\ln R_{n})^{\beta }}%
=n^{a}2^{n\theta }\times \frac{n^{2\alpha }}{2^{2nm\alpha }(2nm\ln 2-2\ln
n)^{\beta }}\leq \frac{C}{n^{\beta -a-2\alpha }}=\frac{C}{n^{2}}
\end{equation*}%
so that%
\begin{equation*}
\sum_{n\geq n_{\ast }}n^{a}2^{n\theta }\left\Vert y-x_{R_{n}}\right\Vert
_{X}\leq \sum_{n\geq n_{\ast }}\frac{1}{n^{2}}<\infty .
\end{equation*}%
$\square $

\smallskip

We come now back to the ``balance'' discussed in Section \ref{sect-results}.
We recall that for $a>1$ we have denoted $L_{a}(R)=R(\ln R)^{a}.$ We also
considered a general Young function $e\in \mathcal{E}$ and we gave the
hypothesis $H_{q}(k,m,e)$. We give here the statement of this hypothesis in
our abstract setting. Let $\theta \geq 0$ and $y\in Y.$

\smallskip

\textbf{Hypothesis $H(\theta ,m,e)$.} \emph{For $\theta\in \N$, $m\in \N_*$
and $e\in \mathcal{E}$ there exists $a>1$ such that%
\begin{equation*}
\overline{\lim }_{R\rightarrow \infty }\frac{L_{a}(R)^{1+\frac{\theta }{2m}}%
}{R}\beta_e (L_{a}(R)^{d/2m})d_{Y}(y,B_{X}(R))<\infty
\end{equation*}%
}

\medskip

For $\theta =k+q,X=W^{2m+q,2m,e}$ and $Y=W_{\ast }^{k,\infty }$ we obtain $%
H_{q}(k,m,e).$

\begin{lemma}
\label{BALANCE}If $H(\theta ,m,e)$ holds for $y\in Y$ then one may find a
sequence $x_{n}$ such that%
\begin{equation*}
\sum_{n=1}^{\infty }2^{n\theta }\beta _{e}(2^{nd})\left\Vert
y-x_{n}\right\Vert _{Y}+\frac{1}{2^{2nm}}\left\Vert x_{n}\right\Vert
_{X}<\infty .
\end{equation*}
\end{lemma}

\textbf{Proof}. For a suitable $a>1$ and large $R$ we have%
\begin{equation*}
d_{Y}(\mu ,B_{X}(R))\leq \frac{CR}{L_{a}(R)^{1+\frac{\theta }{2m}}\beta
_{e}(L_{a}(R)^{d/2m})}.
\end{equation*}%
We choose $R_{n}=n^{-a}2^{2nm}$ and we take a sequence $x_{n}\in
B_{X}(R_{n}) $ such that
\begin{equation*}
\left\Vert y-x_{n}\right\Vert _{Y}\leq \frac{CR_{n}}{L_{a}(R_{n})^{1+\frac{%
\theta }{2m}}\beta _{e}(L_{a}(R_{n})^{d/2m})}.
\end{equation*}%
By the very definition of $B_{X}(R_{n})$ we have
\begin{equation*}
\sum_{n=0}^{\infty }\frac{1}{2^{2nm}}\left\Vert x_{n}\right\Vert _{X}\leq
\sum_{n=0}^{\infty }\frac{1}{2^{2nm}}R_{n}\leq \sum_{n=0}^{\infty }\frac{1}{%
n^{a}}<\infty .
\end{equation*}%
One also has%
\begin{equation*}
(2m)^{a}\times 2^{2nm}\geq L_{a}(R_{n})=\frac{1}{n^{a}}2^{2nm}(2nm\ln 2-a\ln
n)^{a}\geq 2^{2nm}
\end{equation*}%
the last inequality being true for sufficiently large $n$ (we need $2nm\ln
2-a\ln n\geq n).$ It follows that $L_{a}(R_{n})^{d/2m}\geq 2^{nd}$ and this
yields $\beta _{e}(L_{a}(R_{n})^{d/2m})\geq \beta _{e}(2^{nd}).$ We conclude
that
\begin{align*}
2^{n\theta }\beta _{e}(2^{nd})\left\Vert y-x_{n}\right\Vert _{Y}\leq &
2^{n\theta }\beta _{e}(2^{nd})\times \frac{C}{L_{a}(R_{n})^{\frac{\theta }{2m%
}}\beta (L_{a}(R_{n})^{d/2m})}\times \frac{R_{n}}{L_{a}(R_{n})} \\
\leq & (2m)^{a}\times \frac{C}{(\ln R_{n})^{a}}\leq \frac{C(2m)^{a}}{n^{a}}
\end{align*}%
which shows that the first series is also convergent. $\square $

\end{document}